\documentclass{amsart}
\usepackage{longtable}
\usepackage{booktabs}
\usepackage{ifthen}
\usepackage{amssymb}
\usepackage{euscript}
\usepackage[all]{xy}
\usepackage{varioref}
\usepackage{graphicx}
\usepackage{verbatim}
\usepackage{lscape}
\usepackage{epsfig}
\usepackage{xspace}
\usepackage{boxedminipage}
\usepackage{ifthen}
\newboolean{pdfOutput}
\setboolean{pdfOutput}{false}
\ifthenelse{\boolean{pdfOutput}}
{
}
{
\usepackage{graphicx}
}
\newcommand{\drawfile}[1]
{
  \ifthenelse{\boolean{pdfOutput}}
  {
	\epsfig{figure={#1.pdf},width=2in}
  }
  {
  \epsfig{file={#1.eps},width=2in}%
  }
}
\usepackage{rotating}
\newtheorem{theorem}{Theorem}
\newtheorem{corollary}{Corollary}
\newtheorem{lemma}{Lemma}
\newtheorem{proposition}{Proposition}
\theoremstyle{remark}

\newtheorem{definition}{Definition}
\newtheorem{example}{Example}
\newtheorem{remark}{Remark}

\newcounter{remarkCounter}
\newcommand{\head}[1]{\textbf{#1}}

\newlength{\setBracketHeight}

\newcommand{\LieDer}{\ensuremath{\EuScript L}}
\newcommand{\hook}{\ensuremath{\mathbin{ \hbox{\vrule height1.4pt
        width4pt depth-1pt \vrule height4pt width0.4pt depth-1pt}}}}

\newcommand{\pd}[2]{\ensuremath{\frac{\partial{#1}}{\partial{#2}}}}

\newcommand{\C}[1]{\ensuremath{\mathbb{C}^{#1}}}

\newcommand{\Gr}[2]{\ensuremath{\operatorname{Gr}\left({#1},{#2}\right)}}

\newcommand{\GL}[1]{\operatorname{GL}\left({#1}\right)}

\newcommand{\SL}[1]{\operatorname{SL}\left({#1}\right)}
\newcommand{\PSL}[1]{\mathbb{P}\SL{#1}}
\newcommand{\PGL}[1]{\mathbb{P}\GL{#1}}
\newcommand{\slLie}[1]{\mathfrak{sl}\left({#1}\right)}

\newcommand{\Iso}[2]{ \operatorname{Iso} \left({#1},{#2}\right)}

\newcommand{\Sym}[2]{\ensuremath{\operatorname{Sym}^{#1}\left({#2}\right)}}
\newcommand{\Lm}[2]{\ensuremath{\Lambda^{#1} \left ( {#2} \right )}}
\newcommand{\nForms}[2]{\ensuremath{\Omega^{#1} \left ( {#2} \right
    )}}
\newcommand{\Cohom}[2]{\ensuremath{ H^{#1} \left ( {#2} \right
    )}}

\DeclareMathOperator{\Ad}{Ad}

\newcommand{\Norm}[2]{\nu_{#2}{#1}}

\newcommand{\Proj}[1]{\mathbb{P}^{#1}}

\newcommand{\OO}[1]{
  \ensuremath{
    \mathcal{O}
    \ifthenelse{\equal{#1}{0}}
      {}
      {\left({#1}\right)}
  }
}
\newcommand{\OOp}[2]{
  \ensuremath{
    \mathcal{O}
    \ifthenelse{\equal{#1}{0}}
      {}
      {\left({#1}\right)}
    \ifthenelse{\equal{#2}{1}}
      {}
      {^{\oplus{#2}}}
  }
}

\newcommand{\Gtot}{\ensuremath{G}}
\newcommand{\gtot}{\mathfrak{g}}
\newcommand{\Xstruc}[1]{\ensuremath{{#1}_{0}}}
\newcommand{\Point}{\ensuremath{\text{pt}}} 
\newcommand{\Xpt}[1]{\ensuremath{{#1}^{{\Point}}}}
\newcommand{\Xline}[1]{\ensuremath{{#1}^{\text{line}}}}
\newcommand{\Xptline}[1]{\ensuremath{{#1}^{\text{pt,line}}}}
\newcommand{\Gstruc}{\Xstruc{\Gtot}}
\newcommand{\gstruc}{\Xstruc{\gtot}}
\newcommand{\Gpt}{\Xpt{\Gtot}}
\newcommand{\gpt}{\Xpt{\gtot}}
\newcommand{\Gline}{\Xline{\Gtot}}
\newcommand{\gline}{\Xline{\gtot}}
\newcommand{\Gptline}{\Xptline{\Gtot}}
\newcommand{\gptline}{\Xptline{\gtot}}

\newcommand{\Flag}[1]{\ensuremath{\mathbb{F}\left({#1}\right)}}
\newcommand{\TT}{{}^{\sharp} t}
\newcommand{\NA}{\ast}
\newcommand{\Dual}{\Lambda}

\begin{document}
\title{Rational curves and ordinary differential equations}
\author{Benjamin McKay}
\address{University College Cork \\ National University of Ireland}
\date{\today} 
\thanks{Thanks to Valerii Dryuma, Maciej Dunajski, Mark Fels, and Jo{\"el} Merker for pointing
out references to the literature, particularly Fels's thesis.}
\begin{abstract}
Complex analytic 2nd order ODE systems whose 
solutions close up to become rational curves 
are characterized by the vanishing of an
explicit differential invariant, and 
form an infinite dimensional family of integrable
systems.
\end{abstract}
\maketitle
\tableofcontents
\section{Introduction}
Henceforth all
manifolds and Lie groups are complex, and
all maps, vector bundles, functions, sections of vector bundles,
path geometries, connections and differential equations are holomorphic.
\subsection{The problem}
Any 2nd order system of ordinary differential
equation determines a \emph{path geometry} 
(defined below). Moreover the path
geometry determines the system, and all
path geometries are locally obtained
from 2nd order ODE systems. But the
concept of path geometry is geometric and
independent of coordinates. Therefore
we can consider a path geometry on a manifold
as a global generalization of a
2nd order ODE system. Every 2nd order
ODE system has solutions, and so every path geometry has a family
of curves, called the \emph{integral curves}
of the path geometry.
\begin{definition}
A path geometry is \emph{straight} if it
is locally isomorphic to a path geometry
whose integral curves are rational.
\end{definition}
Our problem: to characterize straight path geometries.
The solution is an explicit local condition,
easy to check. 
\begin{example} The fundamental example
which guides our work is the differential equation
\[
\frac{d^2 y}{dx^2} = 0,
\]
whose solutions are straight lines.
This equation is invariant under translations in both
$x$ and $y$, so that we can quotient to define the equation
on a complex torus. However, it is straight because
it is locally isomorphic to the equation of 
projective lines in projective space.
\end{example}
Our problem of characterizing straight path geometries
is similar to Painlev{\'e}'s problem
on systems with fixed singular points, but the answer is quite different.
\subsection{The solution}
\begin{definition}
For a system of $n \ge 1$ second order ordinary differential
equations
\[
\frac{d^2 y^I}{dx^2} = f^I\left(x,y,\frac{dy}{dx}\right).
\]
in complex variables $x,y^1,\dots,y^n$, and for any
function $g\left(x,y,\dot{y}\right)$, define $dg/dx$ to mean
\[
\frac{dg}{dx} = \pd{g}{x} + \pd{g}{y^I} \dot{y}^I + \pd{g}{\dot{y}^I} f^I\left(x,y,\dot{y}\right).
\]
Define the \emph{Fels torsion}
\cite{Fels:1995} of the system to be:
\[
\Phi^I_J = \phi^I_J - \frac{1}{n} \phi^K_K \delta^I_J
\]
where
\[
\phi^I_J =
\frac{1}{2} \frac{d}{dx} \pd{f^I}{\dot{y}^J}
- \pd{f^I}{y^J}
- \frac{1}{4} \pd{f^I}{\dot{y}^K} \pd{f^K}{\dot{y}^J},
\]
with $f^I=f^I\left(x,y,\dot{y}\right)$.
For a single second order ordinary differential equation (i.e. $n=1$),
say
\[
\frac{d^2 y}{dx^2} = f\left(x,y,\frac{dy}{dx}\right),
\]
clearly the Fels torsion vanishes by definition.
Define the \emph{Tresse torsion} (see \cite{Arnold:1988,Cartan:1938,Tresse:1896}:
\[
\frac{d^2}{dx^2} \pd{^2 f}{\dot{y}^2}
-
4 \frac{d}{dx} \pd{^2 f}{y \partial \dot{y}}
+
\pd{f}{\dot{y}}
\left(
  4 \pd{^2 f}{y \partial \dot{y}}
  -
  \frac{d}{dx} \pd{^2 f}{\dot{y}^2}
\right)
- 3 \pd{f}{y} \pd{^2 f}{\dot{y}^2}
+ 6 \pd{^2 f}{y^2}.
\]
\end{definition}
The Fels torsion depends only on second derivatives
of the functions $f^I$, while the Tresse torsion
depends on derivatives of fourth order.
\begin{theorem}\label{thm:Main}
A path geometry is torsion-free (i.e. the Tresse--Fels 
torsion vanishes) just when it is straight.
\end{theorem}
\section{Examples}
\begin{example}
Lets return to our fundamental example,
\[
\frac{d^2 y}{dx^2} = 0,
\]
whose integral curves are straight lines. Lines analytically
continue to become projective lines in projective space.
The differential equation looks the same throughout
projective space, although one has to make projective
linear changes of variable to see what happens out
at the hyperplane at infinity. There is no global
choice of variable $x$ over which integral curves
can be graphed. Naturally the Tresse--Fels torsion vanishes.
\end{example}
\begin{example}\label{example:Linear}
The straight linear second order systems
with constant coefficients are precisely those of the form
\[
\frac{d^2 y}{dx^2} = A \frac{dy}{dx} + \left( a - \frac{1}{4} A^2 \right)y
\]
where $A$ is any constant complex $n \times n$ matrix, and $a$ any complex scalar.
\end{example}
\begin{example}
None of Painlev{\'e}'s equations are straight.
\end{example}
\begin{example}
Any 2nd order ODE with one dimensional symmetry group
can be brought by coordinate transformation to the form
\[
\frac{d^2 y}{dx^2} = f\left(y,\frac{dy}{dx}\right);
\]
(for proof see Lie \cite{Lie:1883}).
The conditions on $f\left(y,\dot{y}\right)$ under
which this equation is straight form the fourth order
equation
\begin{align*}
0 =&
\dot{y}^2 \pd{^4 f}{y^2 \partial \dot{y}^2}
+ 2 \dot{y} f \pd{^4 f}{y \partial \dot{y}^3}
+f^2 \pd{^4 f}{\dot{y}^4}
+ \dot{y} \pd{^3 f}{\dot{y}^3} \pd{f}{y}
- 3 \pd{^3 f}{y \partial \dot{y}^2}
-4 \dot{y} \pd{^3 f}{y^2 \partial \dot{y}}
+ 4 \pd{f}{\dot{y}} \pd{^2 f}{y \partial \dot{y}} \\
&- \dot{y} \pd{f}{\dot{y}} \pd{^3 f}{y \partial \dot{y}^2}
- 3 \pd{f}{y} \pd{^2 f}{y \partial \dot{y}}
+ 6 \pd{f}{y^2},
\end{align*}
so that the generic straight equation with one dimensional 
symmetry group depends
on 4 functions of 1 variable.
In particular, the generic equation with one dimensional symmetry 
group is not straight,
and vice versa; straightness is independent of symmetry group.
\end{example}
\begin{example}
Straightness for 2nd order ODE systems is also independent of
linearizability (see Merker \cite{Merker:2004}).
Any linear 2nd order equation $\frac{d^2 y}{dx^2} = a(x) \frac{dy}{dx} + b(x)y$
is straight. However, the generic coupled system of linear 2nd order equations is
not straight. Consider a single oscillator
\[
\frac{d^2 y}{dx^2} = \omega^2 y,
\]
with $\omega$ constant. The integral curves are
\[
y = a_+ e^{\omega x} + a_- e^{- \omega x},
\]
so that if we introduce the variable $X = e^{\omega x}$, then
\[
y = a_+ X + \frac{a_-}{X},
\]
or
\[
a_+ X^2 - Xy + a_- = 0,
\]
quadratic equations, so the (generic) solutions are smooth rational curves.
Consider a coupled system, say
\begin{align*}
\frac{d^2 y_1}{dx^2} &= \omega_1^2 y_1 \\
\frac{d^2 y_2}{dx^2} &= \omega_2^2 y_2.
\end{align*}
If the frequencies $\omega_j$ are rational multiples
of a common frequency $\omega$, then we can introduce
a parameter $X = e^{\omega x}$, and obtain algebraic equations
for the solutions. But the degrees are not low enough
to keep the curves rational, unless the frequencies $\omega_j$
are all equal: a system of uncoupled
harmonic oscillators is straight just when all
of the frequencies are equal. In consonance with
our theorem (and with example~\vref{example:Linear}), the Fels torsion
vanishes just for equal frequencies.
\end{example}
\begin{example}
For example consider the equation
\[
\frac{d^2 y}{dx^2} = 6 y^2,
\]
which is not torsion-free. Check that the
function $\dot{y}^2 - 4y^3$
is constant along integral curves.
Therefore the integral curves are precisely
the curves
\[
\dot{y}^2 = 4 y^3 + A
\]
for any constant $A$. These curves are elliptic
curves (hence not rational), filling out the phase space, except for the
curve with $A=0$, which is a cuspidal cubic
curve, hence rational; see figure~\vref{fig:cubics}.
\begin{figure}
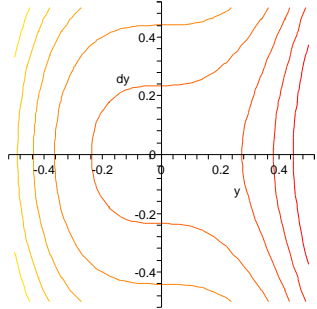

  \centering\drawfile{ellipCurves}
\caption{The family of cubic curves in the plane}%
\label{fig:cubics}
\end{figure}
Integrating, we find
\[
\int \frac{dy}{\sqrt{4 y^3+A}} = x + B,
\]
an elliptic integral on each elliptic curve; the
constant $B$ just translates the elliptic
curve along the $x$ variable, and
is defined up to periods. Going backwards,
$y(x)$ is an elliptic function on each elliptic curve;
in fact it is the Weierstra\ss{}\xspace
$\wp$-function:
$y(x)=\wp\left(x-c\right)$ with modular parameters $g_2=0$ and $g_3=A$.
Globally we can analytically continue all of the integral curves
to the 3-manifold $\C{} \times \Proj{2}$, with
$x$ a coordinate function on $\C{}$, and $\left(y,\dot{y}\right)$
an affine chart on $\Proj{2}$ and each  integral
curve is an elliptic curve, except for the 1-parameter
family of curves with $A=0$ and $B$ arbitrary. The picture
extends to the line at infinity on $\Proj{2}$,
because the elliptic curves are smooth there too.
We have to avoid the surface of points
$\left(x,y,\dot{y}\right)$ for $\left(y,\dot{y}\right)$
in the cuspidal cubic, where the curves don't behave nicely.
\end{example}

\section{Path geometry}
A \emph{path geometry} on a manifold usually
means a differential
system locally given by a 2nd order ODE system,
so that
through any point, in each direction, there is a unique
immersed curve (called an \emph{integral curve})
solving the ODEs passing through that point tangent
to that direction.
In local coordinates $x,y^1,\dots,y^n$ on the manifold,
the integral curves are the solutions of an
equation
\[
\frac{d^2y}{dx^2} = f\left(x,y,\frac{dy}{dx}\right).
\]
Hitchin \cite{Hitchin:1982} shows that
complex surfaces containing rational
curves provide a source of path geometries.
He demonstrates that straight path geometries play a
role in the Penrose twistor
programme.

We need a slightly broader definition of
path geometry:
\begin{definition}\label{def:pathGeometry}
A \emph{path geometry} on a complex manifold $M^{2n+1}$
is a foliation $I$ by curves and a transverse foliation $S$ by $n$-folds,
whose leaves are  respectively called the \emph{integral curves} and \emph{stalks},
so that, if we think of $S$ and $I$ as subsheaves of the tangent 
bundle, $I+S+[I,S]=TM.$
\end{definition}
Near any point of $M$, there is a 
coordinate chart with coordinate functions 
$x,y_1,\dots y_n, \dot{y}_1, \dots, \dot{y}_n$
and there are functions $f^I(x,y,\dot{y})$
in which the integral curves intersect the coordinate chart 
precisely in the solutions of
\[
dy^I = \dot{y}^I \, dx, \ d\dot{y}^I = f^I\left(x,y,\dot{y}\right)  \, dx,
\qquad I=1, \dots, n,
\]
and the stalks intersect the coordinate chart precisely in 
the solutions of
\[
dy^I = dx = 0,
\qquad I=1, \dots, n.
\]

We shall refer to the space of pointed lines in projective
space as the \emph{model}. Its integral curves
are the curves given by moving a point along a fixed line.
Its stalks are given by fixing a point,
and looking at all lines through it.

For any path geometry, following Cartan's terminology, 
call the space
of stalks the \emph{space of points}, and our original
manifold $M$ the \emph{space of elements}.
Locally, every path geometry has
a smooth space of points, with coordinates $x,y$,
but I shall not require the space of points to be smooth globally.
Even if the space of points is smooth,
the path geometry may appear on it as a multivalued
ordinary differential equation, in local coordinates,
and there might be no paths in certain directions
(i.e. $f(x,y,\dot{y})$ might not be defined
for certain values of $\dot{y}$).
I shall mollify this multivaluedness only slightly by
assuming that the space of elements is connected.
I shall not require existence or uniqueness of an 
integral curve in each
direction at a given point $x,y$ in the space of points.

One motivation for this paper is that (as we shall see)
both stalks and integral curves are
canonically locally identified with projective spaces,
modulo projective transformations.
This might remind us of Riemannian geometry, where
geodesics are canonically equipped with arclength
parameterization, defined up to choice
of a constant; the Riemannian manifold is complete
just when the parameterization is a covering of the geodesic by
the real line. The geometry of more general 2nd
order equations is more complicated, and more slippery,
so we have local projective parameterizations defined
only up to projective transformation.  
We shall say integral curves or stalks are
\emph{rational} if they are covered by
projective spaces. For integral curves, this is the natural analogue
of completeness.

We shall prove:
\begin{theorem}\label{theorem:RatStalk}
A stalk [integral curve] of a path geometry
is rational just when it is compact with finite fundamental
group. Moreover this occurs just when the canonical
local identifications with projective spaces extend globally
to a diffeomorphism.
\end{theorem}

Therefore rationality (of the leaves of either foliation)
is a topological condition, but with
strong global consequences.
We shall prove:
\begin{theorem}\label{thm:Model}
The only path geometry on any
connected complex manifold whose integral
curves and stalks are all rational is the path geometry on 
projective space whose integral curves are projective lines.
\end{theorem}

Summing up, we have a topological
criterion for isomorphism with the model.
Note that we do not assume that our complex
manifold is compact or K{\"a}hler.
We shall also prove:
\begin{theorem}\label{theorem:Local}
A path geometry is locally isomorphic
to the path geometry of $d^2 y/dx^2=0$
just when it is both (1) locally isomorphic
to a path geometry with rational integral
curves, and (2) also locally isomorphic to
a path geometry with rational stalks.
\end{theorem}
We shall also locally characterize path geometries
with rational integral curves and those with
rational stalks. Some well known results remain valid
even with our broader definition of
path geometry.
\begin{theorem}[LeBrun \cite{Hitchin:1982}]
If all of the stalks are rational, then (1) 
there is a smooth space of points which
bears a projective connection, (2) the space
of elements is invariantly mapped by local biholomorphism to the
projectived tangent bundle of the space
of points, and (3) each integral curve
is locally identified with the family of tangent lines to a unique
geodesic of the projective connection.
Conversely, every projective connection
on any manifold gives rise to a
path geometry on its projectivized tangent bundle,
with rational stalks (the projectived tangent
spaces). A path geometry is locally
isomorphic to a path geometry with rational
stalks, and therefore is locally a projective
connection, just when it satisfies
\[
\pd{^4 f^i}{\dot{y}^I \dot{y}^J \dot{y}^K \dot{y}^L} = 0,
\]
for any four indices $I,J,K,L=1,\dots,n$, i.e. has the form
\[
\frac{d^2 y}{dx^2} = \sum_{|\alpha| \le 3}
f_{\alpha} \left(x,y\right) \left(\frac{dy}{dt}\right)^{\alpha},
\]
with $\alpha$ a multi-index $\alpha=\left(\alpha_1,\dots,\alpha_n\right).$
\end{theorem}

Hitchin \cite{Hitchin:1982}, Bryant, Griffiths \& Hsu \cite{Bryant/Griffiths/Hsu:1995},
Fels \cite{Fels:1995}, and Grossman \cite{Grossman:2000},
use a more restrictive definition of path geometry,
requiring that there be a smooth space of points
and a smooth space of integral curves; we do not require
either of these, but the reader can easily see that those authors
did not employ these hypotheses in their calculations, only in their conclusions.
\section{Elementary remarks on linearization}
Recall the concept of linearization of a system
of ordinary differential equations: given a
system
\[
\frac{d^2 y^I}{dx^2} = f^I\left(x,y,\frac{dy}{dx}\right),
\]
we linearize about a point $\left(x,y,\dot{y}\right)$ by first
taking the solution $y=y(x)$ through that point,
and then changing coordinates so that the solution becomes
just $y(x)=0$, and the point becomes $(0,0,0)$,
and then we expand $f^I$ into a Taylor
expansion and keep the lowest order terms. It is
thus elementary to see that
\begin{theorem}
A system of 2nd order ordinary differential
equations is torsion-free just when its linearization
about any point is torsion-free, which occurs just when its linearization
about any point has the form
\[
\frac{d^2 y}{dx^2} = A \frac{dy}{dx} + \left( a - \frac{1}{4} A^2 \right)y.
\]
\end{theorem}
\section{A first glance at surface path geometries}
A path geometry shall be called
a \emph{surface} path geometry to denote
that there is one $y$ variable (and there is always
only one $x$ variable), i.e. that the space of points is a (not
necessarily Hausdorff) surface. Therefore the
space of elements $M$ is a 3-fold. 
\begin{theorem}
If the stalks [integral curves] of a surface path geometry are compact,
and there is a 
submersion from a nonempty open set $\pi : U \text{ open } \subset M^3 \to S^2$ and
whose fibers are stalks [integral curves] respectively, and
the space of elements $M$ is connected, then all of the stalks [integral curves]
are rational.
\end{theorem}
\begin{proof}
Take each
point $m \in U$, construct the integral
curve $C_m \subset M$ through $m$,
and map $m \in U \mapsto \pi'(m) T_m C_m \in \Proj{}TS$.
In local coordinates $x,y,\dot{y}$, evidently
this is a local biholomorphism, mapping stalks
to fibers $\Proj{}T_s S$. By compactness of stalks,
this map is onto. Stalks are connected by definition,
so the map is a covering map on each stalk, and therefore
a biholomorphism, because the fibers $\Proj{}T_s S = \Proj{1}$
are simply connected. By analytic continuation, 
all stalks are rational.
\end{proof}
\subsection{The structure equations}
I draw freely from Bryant, Griffiths \& Hsu \cite{Bryant/Griffiths/Hsu:1995}.
They prove that given any surface path geometry,
$M$ bears a canonical choice of 
Cartan geometry, which we write as 
$E \to M$ with Cartan connection $\omega$, modelled 
on the space of points lines in the projective plane.
(They give $E$ the name $B_{G_3}$).
So $E$ is a principal right $\Gptline$-bundle,
where $\Gptline \subset \Gtot=\PGL{3,\C{}}$
is the subgroup fixing a projective line in the projective plane and
a point on that line, i.e. the
group of matrices of the form
\[
g = \begin{bmatrix}
g^0_0 & g^0_1 & g^0_2 \\
   0 & g^1_1 & g^1_2 \\
   0 & 0     & g^2_2
\end{bmatrix},
\]
with Lie algebra $\gptline$. (The square brackets indicate that the matrix
is defined up to rescaling, being an element
of $\PGL{3,\C{}}$). Moreover, they define
a canonical 1-form $\omega$ (which they write
as $\phi$) on $E$ valued in
$\gtot=\slLie{3,\C{}}$ (the Lie algebra of $\Gtot=\PGL{3,\C{}}$),
so that
\begin{enumerate}
\item
$\omega_e : T_e E \to \gtot$ is a linear
isomorphism
\item\label{eqn:BGHstruc}
$\omega \pmod{\gptline}$ is semibasic
for $E \to M$, and
\[
d \omega = - \frac{1}{2} \left[\omega,\omega\right] +
\nabla \omega
\]
where (writing $\omega = \left(\omega^i_j\right)$)
\[
\nabla \omega
=
\begin{pmatrix}
0 & K_1 \omega^1_0 \wedge \omega^2_0 & \omega^2_0 \wedge \left( L_1 \omega^1_0 + L_2 \omega^2_1 \right) \\
0 & 0 & K_2 \omega^2_0 \wedge \omega^2_1 \\
0 & 0 & 0
\end{pmatrix}.
\]
Write $r_g : E \to E$ for the right action of
an element $g \in \Gptline$ on $E$.
\item\label{item:BGH}
\[
r_g^* \omega = \Ad_g^{-1} \omega,
\]
\item\label{item:Lift}
Given any local section $\sigma$ of $E \to M$, the
integral curves are precisely the
solutions of the exterior differential system $\sigma^* \omega^2_0=\sigma^* \omega^2_1=0$.
\end{enumerate}
Bryant, Griffiths \& Hsu don't actually state the equation
(\ref{item:BGH}), but it
follows immediately as a simple calculation
from the transformation properties of the various components of $\omega$
as given in their article. They also don't state
(\ref{item:Lift}), but it is clear from their remarks
on the top of p. A.2.
\section{Review of Cartan connections}
In this section, I define Cartan connections,
and prove a few results
which were only sketched in my paper \cite{McKay:2004}.
\begin{definition}
  A \emph{Cartan pseudogeometry} on a manifold $M$, modelled on a
  homogeneous space $\Gtot/\Gstruc$, is a principal right $\Gstruc$-bundle $E \to M$,
  (with right $\Gstruc$ action written $r_g : E \to E$ for $g \in \Gstruc$),
  with a 1-form $\omega \in \nForms{1}{E} \otimes \gtot$,
  called the \emph{Cartan pseudoconnection} (where $\gtot, \gstruc$ are
  the Lie algebras of $\Gtot,\Gstruc$), so that $\omega$ identifies each tangent
  space of $E$ with $\gtot$.  For each $A \in \gtot$,
  let $\vec{A}$ be the vector field on $E$ satisfying $\vec{A} \hook
  \omega = A$.  A Cartan pseudogeometry is called a \emph{Cartan
    geometry} (and its Cartan pseudoconnection called a \emph{Cartan
    connection}) if (1) $r_g^* \omega = \Ad_g^{-1} \omega$ for all $g \in \Gstruc$
    and (2)
    \[
    \vec{A} = \left. \frac{d}{dt}r_{e^{tA}}\right|_{t=0}
    \]
    for all $A \in \gstruc.$
\end{definition}
\begin{lemma}
The 1-form $\omega$ of Bryant, Griffiths \& Hsu is a Cartan connection on $M$,
modelled on $\Gtot/\Gptline=\Proj{}T\Proj{2}=\Flag{1,2}$, the flag variety
of pointed lines in projective space.
\end{lemma}
\begin{proof}
We have only to check that $\vec{A}$ is the infinitesimal
generator of the right action, for $A \in \gptline$.
This follows immediately from the simple calculation
that
\[
\LieDer_{\vec{A}} \omega = - \left[A,\omega\right].
\]
\end{proof}
We employ a host of results on vector bundles and Cartan geometries,
all of which have the same proof, so we give the
proof in just one case:
\begin{lemma}\label{lemma:TgtBundle}
Consider a Cartan geometry $\pi : E \to M$. The tangent bundle is
\[
TM = E \times_{\Gstruc} \left( \gtot/\gstruc \right ).
\]
\end{lemma}
\begin{proof}
At each point $e \in E$, the 1-form $\omega_e : T_e E \to \gtot$
is a linear isomorphism, taking $\ker \pi'(e) \to \gstruc$. Therefore
$\omega_e : T_e E / \ker \pi'(e) \to \gtot/\gstruc$ is a linear
isomorphism.
Also $\pi'(e) : T_e E /\ker \pi'(e) \to T_{\pi(e)} M$
is an isomorphism. Given a function $f : E \to \gtot/\gstruc$,
define $v_f$ a section of the vector bundle $TE / \ker \pi'$,
by the first isomorphism, and a section $\bar{v}_f$
of $\pi^* TM$ by the second. Calculate that $\bar{v}_f$
is $\Gstruc$-invariant just when $f$ is $\Gstruc$-equivariant,
i.e. just when
\[
r_g^* f = \Ad_g^{-1} f.
\]
This makes an isomorphism of sheaves between the
sections of the tangent
bundle $TM$ and the $\Gstruc$-equivariant functions
$E \to \gtot/\gstruc$, i.e. the sections of
$E \times_{\Gstruc} \left(\gtot/\gstruc\right)$,
so that they must be identical vector bundles.
\end{proof}
\begin{definition}
For $G_0 \subset G$ a closed subgroup, let
$\omega \in \nForms{1}{G}$ be the left invariant
Maurer--Cartan 1-form. Then $\omega$ is a Cartan
connection on the principal right $G_0$-bundle
$G \to G/G_0$, and the induced Cartan geometry
on $G/G_0$ is called the model Cartan geometry.
A Cartan geometry modelled
on $G/G_0$ is called \emph{flat} if it is locally
isomorphic to the model Cartan geometry.
\end{definition}
\begin{definition}
The expression $\nabla \omega = d \omega + \frac{1}{2} [ \omega, \omega ]$
is called the \emph{curvature} of the Cartan geometry;
equations on the curvature are called \emph{structure equations}.
\end{definition}
\begin{theorem}[Sharpe \cite{Sharpe:1997}]
A Cartan geometry is flat just when its curvature vanishes.
\end{theorem}
\begin{proposition}
Pick a flat Cartan geometry $E \to M$ on a compact,
connected and simply connected manifold $M$, modelled on
$G/G_0$ with $G$ connected and $G/G_0$ connected
and simply connected. Then the
Cartan geometry is isomorphic to the model.
\end{proposition}
\begin{proof}
By theorem 3 of
McKay \cite{McKay:2004}, some covering space of $M$
maps locally diffeomorphically to $G/G_0$,
and the Cartan geometry
on that covering space is pulled back.
Because $M$ is simply connected, that
covering space is $M$ itself. Because
$M$ is compact, the local diffeomorphism
is a covering map. Because $G/G_0$ is
connected and simply connected, the map
is a diffeomorphism.
\end{proof}
\begin{definition}
If $\Gstruc \subset \Gtot$ is a closed
subgroup of a Lie group, and $\Gamma \subset \Gtot$
is a discrete subgroup, acting freely
and discontinuously on $\Gtot/\Gstruc$,
then we can let
$E=\Gtot, M=\Gamma \backslash \Gtot/\Gstruc, \omega=g^{-1} , dg$,
determining a flat Cartan geometry called
a \emph{locally Klein geometry}.
\end{definition}
Say that a group $G$ \emph{defies} a group $H$
if every morphism $G \to H$ has finite image.
We do not repeat the proof of:
\begin{theorem}[McKay \cite{McKay:2004}]\label{thm:Defiance}
A flat Cartan geometry, modelled on $\Gtot/\Gstruc$, defined on a compact
connected base manifold $M$ with fundamental group defying $G$, is a locally Klein
geometry.
\end{theorem}
\begin{definition}
If $V$ is a vector space, a $V$-valued \emph{coframing} 
on a manifold $E$ is
1-form $\omega \in \nForms{1}{E} \otimes V$,
so that at each point $e \in E$, $\omega_e : T_e E \to V$ 
is a linear isomorphism.
An isomorphism of coframings is a diffeomorphism
matching up the 1-forms. 
If $\Gstruc \subset \Gtot$ is a closed Lie subgroup
of a Lie group, with Lie algebras $\gstruc \subset \gtot$,
and $\omega$ is a $\gtot$-valued coframing,
let 
$\bar{\omega} = 
\omega \mod \gstruc \in \nForms{1}{E} \otimes \left( \gtot/\gstruc \right).$
A \emph{local Cartan geometry} modelled
on $\Gtot/\Gstruc$ on a manifold $E$ is a $\gtot$-valued coframing $\omega$ on $E$
and a function $K : E \to \Lm{2}{\gtot/\gstruc} \otimes \gtot$,
for which
\[
d \omega + \frac{1}{2} \left[\omega,\omega\right] = K \bar{\omega} \wedge \bar{\omega}.
\]
\end{definition}
\begin{definition}
If $E \to M$ bears a Cartan geometry with Cartan connection
$\omega$, then $\omega$ and the curvature $K$ of $\omega$
are together called the \emph{associated local Cartan geometry}.
We say that a local Cartan geometry 
is isomorphic to a Cartan geometry if it is isomorphic to 
the associated local Cartan geometry.
\end{definition}
\begin{theorem}\label{thm:local}
Every local Cartan geometry is locally isomorphic
to a Cartan geometry.
\end{theorem}
\begin{remark}
This theorem is a well-known folk theorem, but we
know of no source for a proof.
\end{remark}
\begin{proof}
Consider the foliation of $E$ by the submanifolds $\bar{\omega}=0$.
Since our result is local, we can assume that this
foliation is a fiber bundle $E \to M$, and also
that this fiber bundle is trivial.
Consider the vectors fields $\vec{A}_E$ on $E$ defined
by the equation $\vec{A}_E \hook \omega = A$,
for any $A \in \gstruc$. These vector fields generate
an action of the Lie algebra $\gstruc,$ whose orbits are the fibers
of $E \to M$. Taking any local section of $E \to M$,
say $\sigma : M \to E$, the map
\[
(m,A) \in M \times \gstruc \mapsto e^A m \in E
\]
is defined near $A=0$, and a local diffeomorphism
there. Therefore we can find an open set of the
form $U_M \times U_{\gstruc}$ with $U_M \subset M$
and $U_{\gstruc} \subset \gstruc$ open sets,
on which the map is defined and is a diffeomorphism to its image.
Because our results are local, we can assume that
$M=U_M$, and the map is a global diffeomorphism, with image
all of $E$. Moreover, we can assume that the
exponential map identifies $U_{\gstruc}$ with
an open subset $U_{\Gstruc} \subset \Gstruc$. We therefore
have $E = M \times U_{\Gstruc} \subset M \times \Gstruc$.
On $M \times \Gstruc$, define
a 1-form $\Omega$ by
\[
\Omega_{\left(m,g_0\right)}=\Ad_{g_0}^{-1} \omega_{(m,1)} r_{g_0^{-1}}'\left(m,g_0\right).
\]
Check that $\Omega=\omega$ on $M \times 1$ and that
$\LieDer_{\vec{A}} \Omega = -[A,\Omega]$, for $A \in \gstruc$,
so that by uniqueness of solutions
of ordinary differential equations, $\Omega=\omega$ on $E$.
The coframing $\Omega$ is a Cartan geometry on $M \times G_0.$
\end{proof}
\section{Inducing a Cartan connection on integral curves}
Let $C \to M$ be any immersed integral curve of a path geometry.
Consider the pullback subbundle $\left. E \right|_C$.
Since $0=\omega^2_0=\omega^2_1$ along $C$ on every local section
of $E \to M$, and $0=\omega^2_0=\omega^2_1$
on the fibers, we find that $0=\omega^2_0=\omega^2_1$
on all of $\left. E \right|_C$.
Moreover, $\left. E \right|_C \to C$ is a principal right $\Gptline$-bundle.
\begin{lemma}
On $\left. E \right|_C \to C$, $\omega$ is a flat Cartan
connection.
\end{lemma}
\begin{proof}
The structure equations are identical to those of $E \to M$, except that $\nabla \omega=0$
because $\omega^2_0=\omega^2_1=0$.
\end{proof}
\section{Classification of Cartan connections on rational curves}
\begin{definition}
A \emph{projective representation} is 
a morphism of complex Lie groups $\alpha : \Gtot \to \PGL{n+1,\C{}}$.
A projective representation is \emph{transitive} if
$\Gtot$ acts transitively on $\Proj{n}$.
Given a transitive projective representation,
set $\Gstruc= \ker \alpha$, $E=\Gtot$,
and $\omega= g^{-1} \, dg$
the left invariant Maurer--Cartan 1-form on $\Gtot$.
Call this the Cartan geometry associated
to the transitive projective representation.
\end{definition}
\begin{theorem}\label{thm:TransProjRep}
Every flat Cartan geometry on $\Proj{n}$,
with connected model $\Gtot/\Gstruc$, is isomorphic to its model, hence
isomorphic to the Cartan geometry associated
to a transitive projective representation.
\end{theorem}
\begin{proof}
Any flat Cartan geometry on $\Proj{n}$
is obtained by taking a local
biholomorphism to the model $\Proj{n} \to \Gtot/\Gstruc$,
so a covering map (since $\Proj{n}$ is compact).
The deck transformations must be biholomorphisms of $\Proj{n}$, so projective
linear transformations. However, every projective
linear transformation has a fixed point, so only
the identity map can act as a deck transformation.
\end{proof}
\begin{corollary}\label{cor:ProjRep}
Every Cartan geometry on a rational
curve is associated to a transitive surjective
projective representation.
\end{corollary}
\begin{proof}
Any Cartan geometry on a curve is flat,
since the curvature is a semibasic 2-form.
No complex Lie subgroup of $\PGL{2,\C{}}$ acts
transitively on $\Proj{1}$. Therefore the projective
representation $\Gtot \to \PGL{2,\C{}}$ is surjective.
\end{proof}
\section{A cornucopia of vector bundles}
If $M^3$ bears a surface path geometry, then
the stalks are curves transverse to the integral
curves. Let $\Theta \subset TM$ be the field of 2-planes
spanned by the tangent lines to integral curves
and tangent lines to stalks. In local coordinates,
$x,y,\dot{y}$, we see that $\Theta=\left(dy=\dot{y} \, dx\right)$,
so a contact structure.
\begin{proposition}
Let $C$ be an immersed integral curve
in a complex 3-fold $M$ with path geometry. 
Let $E \to M$ be the Cartan geometry
associated to the path geometry.
Let $\Theta \subset M$ be the
canonical contact structure.
Let $S$ be the space of points. (If
$S$ is not a smooth surface, then equations
below involving $S$ are meaningless,
but the right hand sides still define vector bundles.)
Let $\Norm{C}{S}=\left.TS\right|_C/TC$ be the normal bundle of the 
immersion $C \to M \to S$
(for which a similar proviso applies).
Let $\Gtot=\PGL{n+1,\C{}}$, $\Gpt$ the subgroup preserving
the point $\left[e_0\right]$, $\Gline$ the subgroup preserving
the line through $\left[e_0\right]$ and $\left[e_1\right]$,
$\Gptline$ the subgroup preserving the point and the line,
and write their Lie algebras as $\gtot, \gpt$, etc.
Then
\begin{align*}
\left. TM \right|_C &= 
\left. E \right|_C \times_{\Gptline} \left(\gtot/\gptline \right) \\
\left. \Theta \right|_C &= 
\left. E \right|_C \times_{\Gptline} 
\left(
	\left(\gline + \gpt\right)
/\gptline 
\right) \\
\Norm{C}{M} &= \left. E \right|_C \times_{\Gptline} \left(\gtot/\gline\right) \\
\left .
	T S
\right |_C &=
    \left.
		E
	\right|_C
	\times_{\Gptline}
	\left(
		\gtot / \gpt
	\right)
\\
TC &= \left. E \right|_C \times_{\Gptline} 
\left(\gline/\gptline \right) \\
\Norm{C}{S} &=
	\left. E
	\right|_C \times_{\Gptline}
	\left(
		\gtot  /\left(\gpt + \gline\right)
	\right) \\
\end{align*}
\end{proposition}
\begin{proof} 
As in the proof of lemma~\ref{lemma:TgtBundle}.
\end{proof}
Intuitively, the equations above allow us to pretend to
work with the space of points $S$, even if it isn't Hausdorff,
by instead working with various vector bundles.

Recall that $\Proj{1}$ bears line bundles $\OO{p}$,
defined as follows: think of $\Proj{1}$ as the space
of lines through 0 in $\C{2}$, and let $\OO{-1}$
be the bundle whose fiber above a line $L$ is
just $L$; then let $\OO{p}=\OOp{-1}{-p}$.
So a local section of $\OO{p}$ is a choice of
map $z \in \text{open } \subset \C{2} \backslash 0 \to f(z)z$
for which $f(\lambda z)=\lambda^p f(z)$.
Moreover, the line bundles $\OO{p}$ have global
nonzero sections just when $p > 0$.
Another way to present these line bundles:
Let
\[
e_0 = \begin{pmatrix}1 \\ 0\end{pmatrix}\in \C{2} \backslash 0.
\]
Write $B$ for the group
of matrices of the form:
\[
g_0 =
\begin{pmatrix}
a & b \\
0 & 1/a
\end{pmatrix}
\]
with $a \ne 0$ (i.e. those matrices 
which preserve the complex line through $e_0$).
Consider the principal right $B$
bundle $\SL{2,\C{}} \to \Proj{1}$ given by the map
$g \in \SL{2,\C{}} \mapsto ge_0 \in \Proj{1}$.
Given an open subset $U \subset \Proj{1}$,
let $\SL{2,\C{}}_U \to U$ be the pullback bundle.
\begin{lemma}\label{lemma:SheafIso}
Sections of $\OO{p}_U \to U$
correspond to maps $F : \SL{2,\C{}}_U \to \C{}$
for which
\[
F\left(gg_0\right) = a^p F(g),
\]
for all $g_0 \in B$.
\end{lemma}
\begin{proof}
Pick a local section $f$ of $\OO{p}$, i.e.
a choice of map $f : \hat{U} \subset \C{2} \backslash 0 \to \C{}$, where
$\hat{U}$ is the preimage of $U$ under $\C{2} \backslash 0 \to \Proj{1}$,
and with $f(az)=a^p z$. Define $F(g) = f\left(ge_0\right)$. Conversely, given
$F$, define $f(z)=F\left(g_z\right)$ where
\[
z=
\begin{pmatrix}
z_1 \\
z_2
\end{pmatrix}, \text{ and }
g_z =
\begin{pmatrix}
z_1 & 0 \\
z_2 & z_1^{-1} \\
\end{pmatrix},
\]
defined for all $z_1 \ne 0$ for which $g_z$ lies
in the domain of $F$.
\end{proof}
In giving this proof, we are merely trying to avoid
abstract Borel--Weil--Bott theory,
and give concrete expressions for these line bundles
``upstairs.''
\begin{corollary}
On any rational integral curve,
\begin{align*}
\left. TM \right|_C &= \OO{2} \oplus \OOp{0}{2} \\
\left. \Theta \right|_C &= \OO{2} \oplus \OO{-1} \\
\left. TM/\Theta \right|_C &= \OO{1} \\
\Norm{C}{M} &= \OOp{0}{2} \\
\left. TS \right|_C &= \OO{2} \oplus \OO{1} \\
TC &= \OO{2} \\
\Norm{C}{S} &= \OO{1} \\
\end{align*}
\end{corollary}
\begin{proof}
Calculate these on the model.
(Note that the normal bundle of $C=\Proj{1} \subset M =\Proj{}T\Proj{2}$
has sections coming from the tangent bundle
of the dual space $\Proj{2*}$; from which it is
easy to see that this normal bundle is trivial.)
By the classification
of Cartan geometries on rational curves, the Cartan
geometry on every rational integral curve is isomorphic to
the one found on the integral curves of the model,
making these vector bundles identical.
\end{proof}
\section{Kodaira deformation theory}
We give a brief review of Kodaira's theory
\cite{Kodaira:1961,Kodaira:1962,Kodaira:2005}.
\begin{definition}
Let $Y$ and $M$ be
complex manifolds and let $\pi_M: M \times Y \to M$ 
and $\pi_Y: M \times Y \to Y$ be the obvious maps.
A \emph{family} of closed complex 
submanifolds of the complex manifold
$M$ parameterized by $Y$ is a complex submanifold 
$F \subset M \times Y$ such that the
$\left.\pi_Y \right|_F : F \to Y$ is a proper submersion.
Let $X_y = F \cap M \times \left\{y\right\}$.
\end{definition}
\begin{definition}
A \emph{morphism} of families $F_j \subset M \times Y_j$
(in the same manifold $M$), $j=0,1$, is a map $\phi : Y_0 \to Y_1$
so that $\left(m,y_0\right) \to \left(m,\phi\left(y_0\right)\right)$
takes $F_0$ to $F_1$.
\end{definition}
\begin{definition}
We say that a submanifold $X \subset M$
belongs to a family $\left\{X_y\right\}_{y \in Y}$ if 
$X=X_y$ for some $y \in Y$.
\end{definition}
\begin{definition}
A family $\left\{X_y \right\}_{y \in Y}$
is \emph{locally complete} if, should one of the
submanifolds $X_y$ belong to another family of
complex submanifolds $\left\{X_z \right\}_{z \in Z}$,
say $X_{y_0}=X_{z_0} \subset M$, then there is morphism 
of families $U \to Y$ defined on an open neighborhood
$U \subset Z$ of $z_0$, taking $z_0 \mapsto y_0$.
\end{definition}
\begin{definition}
A closed complex submanifold $X \subset M$
with normal bundle $\nu_X$ is \emph{free}
if $\Cohom{1}{X,\nu_X}=0$.
\end{definition}
\begin{theorem}[Kodaira \cite{Kodaira:1961,Kodaira:1962,Kodaira:2005}]
If $X \subset M$ is an immersed free closed complex
submanifold of a complex manifold, then $X$
belongs to a locally complete family of submanifolds
$\left\{X_y\right\}_{y \in Y}$,
with an isomorphism $T_X Y=\Cohom{0}{X,\nu_X}$.
If $\Cohom{1}{X,TX}=0$, then every manifold
in this family is biholomorphic to every other.
\end{theorem}
\begin{corollary}
Let $X$ be a closed complex manifold
with $\Cohom{1}{X,TX}=0$.
Let $M$ be a complex manifold, 
and $Y$ the set of free closed complex submanifolds
of $M$ biholomorphic to $X$. Then
$Y$ is either empty or a complex manifold
of dimension equal to the dimension of $\Cohom{0}{X,\nu_X}$,
and a locally complete family.
\end{corollary}
To make use of this, we need to know
a little sheaf cohomology:
\begin{lemma}
\begin{align*}
\dim \Cohom{0}{\Proj{1},\OO{p}} &=
\begin{cases}
p+1 & p \ge 0 \\
0   & p < 0 \\
\end{cases}
\\
\dim \Cohom{1}{\Proj{1},\OO{p}} &=
\begin{cases}
0 & p \ge 0 \\
\left|p+1\right|   & p < 0 \\
\end{cases}
\end{align*}
\end{lemma}
See Griffiths \& Harris \cite{GriffithsHarris:1994} for proof.
\begin{corollary}
If a surface path geometry has
a rational integral curve [rational stalk],
and all integral curves [stalks] compact, then all of its
integral curves [stalks] are rational, and the space
of integral curves [points] is a smooth surface.
\end{corollary}
\section{Identifying line bundles}
Let $E \to C$ be a Cartan geometry on a rational
curve, with Cartan connection $\omega$,
modelled on a homogeneous space $G/G_0$.
Then $E=G,$ and the Cartan connection is the
Maurer--Cartan form, by corollary~\vref{cor:ProjRep}.
Suppose that $\Gstruc$ and $\Gtot$ are connected.
On $\PGL{2,\C{}},$ we can write the Maurer--Cartan form as
\[
g^{-1} \, dg = 
\begin{pmatrix}
\omega^0_0 & \omega^0_1 \\
\omega^1_0 & -\omega^0_0
\end{pmatrix}.
\]
We pull this form back to $E$.
Consider a function $F : \text{connected open } \subset E \to \C{}$.
If there is a number $p$ for which
\[
dF - p F \omega^0_0
\]
is semibasic, then call $p$ the \emph{weight} of $F$.
\begin{lemma}
Suppose that $E \to C$ is a Cartan geometry on a rational curve.
A function $F$ on a connected open subset of $E$
has integer weight $p$ just when $F$ is a local
section of $\OO{p}$ on $C$.
\end{lemma}
The proof is clear from lemma~\vref{lemma:SheafIso}.
\begin{corollary}
Functions $F : E \to \C{}$ of negative integer weight vanish.
\end{corollary}
\section{Dual surface path geometries}
It is an old observation that every surface path
geometry has a dual surface path geometry, given by
interchanging the role of integral curves and stalks;
see 
Bryant, Griffiths and Hsu \cite{Bryant/Griffiths/Hsu:1995}
or
Crampin and Saunders \cite{CrampinSaunders:2005}. We can see
directly from the structure equations of
the Cartan geometry that this is just an interchange
of indices $\omega^{\mu}_{\nu} \leftrightarrow \omega^{2-\nu}_{2-\mu}$.
\begin{example}
For the model, this duality is the duality between
lines in the projective plane and points in the dual
plane, i.e. between 2-planes through 0 in $\C{3}$
and lines through 0 in $\C{3*}$, given by
$\Pi \mapsto \Pi^{\perp}$.
\end{example}
\begin{example}[Hitchin \cite{Hitchin:1982} p. 83]
The ordinary differential equation
\[
\frac{d^2y}{dx^2} = \frac{1}{4y^3}
\]
has solutions
\begin{equation}\label{eqn:Hitchin}
y^2 = ax^2 + bx + c
\end{equation}
for any complex constants $a,b,c$ for which $4ac-b^2=1$.
Now treat $x,y$ as constants, and think of equation~\vref{eqn:Hitchin}
as determining a family of curves $b=b(a),c=c(a)$.
Differentiating twice, we find the relationship:
\[
\frac{d^2 c}{da^2}
=
\frac{4Q(a \frac{dc}{da} + c + \sqrt{Q})}{(4ac-1)(2a^2 \frac{dc}{da} - 2 ac + 2 a \sqrt{Q} + 1)},
\]
the dual ordinary differential equation,
where
\[
Q = c^2 - 2 a c \frac{dc}{da} + a^2 \left( \frac{dc}{da} \right)^2 + \frac{dc}{da}.
\]
The integral curves of the dual equation consist precisely in
the values of $a,b,c$ which produce
a solution $y=y(x)$ of the original equation which passes
through a chosen point of the $(x,y)$ plane.

Even though the equations
\[
y^2 = ax^2 + bx + c, 4ac-b^2=1
\]
are quadratic, so the integral curves and stalks are rational
curves in the plane, the original differential equation
\[
\frac{d^2y}{dx^2} = \frac{1}{4y^3}
\]
is \emph{not} torsion-free, detecting the singularity
emerging in the ordinary differential equation at $y=0$.
Nonetheless, the dual equation is torsion-free,
and is straight.
\end{example}
\begin{example}
The ordinary differential equation
\[
\frac{d^2 y}{dx^2} =
\frac{1}{2} y (y-1)
+
\frac{y-\frac{1}{2}}{y(y-1)}\left(\frac{dy}{dx}\right)^2
\]
conserves the quantity
\[
\lambda = y - \frac{\left(\frac{dy}{dx}\right)^2}{y(y-1)},
\]
from which we conclude that
\[
\frac{dy}{dx} = \sqrt{y(y-1)(y-\lambda)},
\]
i.e. an elliptic curve in phase space, giving
\[
x + \varpi = \int \frac{dy}{\sqrt{y(y-1)(y-\lambda)}},
\]
where $\varpi$ is the integral over a period.
Our elliptic curve has equation
\[
\dot{y}^2 = y(y-1)(y-\lambda).
\]
Compute (as in Clemens \cite{Clemens:2003} p. 59) that
\[
d
\left(
\frac{\dot{y}}{(y-\lambda)^2}
\right)
=
-\frac{1}{2} \frac{dy}{\dot{y}}
- 2(2 \lambda-1) \pd{}{\lambda} \frac{dy}{\dot{y}}
-2 \lambda (\lambda - 1) \pd{^2}{\lambda^2} \frac{dy}{\dot{y}}.
\]
Integrating both sides along the elliptic
curve, avoiding $\dot{y}=0$, we find the Picard--Fuchs equation
\[
0 = \lambda (\lambda-1)\frac{d^2 \varpi}{d \lambda^2}
+ (2 \lambda - 1) \frac{d \varpi}{d \lambda}
+ \frac{1}{4} \varpi,
\]
the dual path geometry, so that $(x,y)$ are the variables on the original
space, and $(\lambda,\varpi)$ are the variables on the dual
space. It is remarkable that the original equation has
nowhere vanishing Tresse torsion (consonant
with our theory, since
all of the integral curves are elliptic curves),
but the Picard--Fuchs equation is torsion-free, so
straight.
\end{example}
\begin{example}
Similarly, returning to our previous example of
\[
\frac{d^2 y}{dx^2} = 6y^2,
\]
the solutions are given implicitly by
\[
\int \frac{dy}{\sqrt{4y^3 + A}} = x + B,
\]
and the dual equation describes how to vary $B$ as a function
of $A=g_3$, in order to keep this elliptic function passing
through a fixed point $x_0,y_0$, i.e. the relation between
the modular parameter $g_3$ (with $g_2=0$)
and the period of the elliptic integral. This is equivalent
to solving for $B$ as a function of $A$ in the equation
\[
y_0 = \left.\wp\left(x_0+B\right)\right|^{g_2=0}_{g_3=A}.
\]
Again there is a Picard--Fuchs equation, but it is more difficult
to find, and we will not try to find it.
\end{example}
\section{Rationality of integral curves and stalks
for surface path geometries}
\begin{theorem}[Hitchin \cite{Hitchin:1982}]
If the integral curves [stalks] of a 
surface path geometry are rational, then $K_1=0$ [$K_2=0$].
\end{theorem}
\begin{proof}
Calculating the exterior derivatives of the structure
equation in item~\vref{eqn:BGHstruc}, we find
\begin{align*}
\nabla
\begin{pmatrix}
K_1 \\
K_2 \\
\end{pmatrix}
&=
d
\begin{pmatrix}
K_1 \\
K_2 \\
\end{pmatrix}
+
\begin{pmatrix}
5 \omega^0_0 + 3 \omega^1_1 & 0 \\
0 & 4\omega^0_0 + 5 \omega^1_1 \\
\end{pmatrix}
\begin{pmatrix}
K_1 \\
K_2 \\
\end{pmatrix}
\\
&=
\begin{pmatrix}
\nabla^0_1 K_1 \omega^1_0 + \nabla^0_2 K_1 \omega^2_0 + L_1 \omega^2_1 \\
\nabla^0_1 K_2 \omega^1_0 + L_2 \omega^2_0 + \nabla^1_2 K_2 \omega^2_1 \\
\end{pmatrix}.
\end{align*}
If $C$ is a rational integral
curve, the bundle $\left.E\right|_C \to C$ has
$\omega^2_0=\omega^2_1=0$. Consider
the copy of $\slLie{2,\C{}}$
determined by
\[
0=\omega^{\bullet}_2=\omega^2_{\bullet}=\omega^0_0+\omega^1_1.
\]
Our equations simplify to
\[
d
\begin{pmatrix}
K_1 \\
K_2
\end{pmatrix}
=
\begin{pmatrix}
-2 & 0 \\
0 & 1
\end{pmatrix}
\omega^0_0
\begin{pmatrix}
K_1 \\
K_2 \\
\end{pmatrix}
\]
modulo semibasic terms, hence weights
$-2,1$ for $K_1, K_2$ respectively.
Therefore $K_1=0$. Duality proves the results for stalks.
\end{proof}
\begin{corollary}
If the space of elements of a surface path
geometry is connected, and both the generic integral curve
and the generic stalk are rational, then the path geometry is
locally isomorphic to the model.
\end{corollary}
\begin{proof}
We see that $K_1=K_2=0$, and differentiate to find
that all invariants vanish. 
The structure equations are identical
to those of the model. Let
$\omega$ be the Cartan connection of the
path geometry on $E \to M$, and $\omega'$ the Cartan connection
of the model, say $E' \to M'$. Then on $E \times E'$,
the holonomic differential system $\omega=\omega'$
has as integral manifolds the graphs of local
isomorphisms.
\end{proof}
\begin{proposition}[Cartan \cite{Cartan:1938}]
Under any local choice of section of $E \to M$,
$K_1$ is Tresse torsion, up to multiplication
by a nowhere vanishing function. In particular,
Tresse torsion vanishes just when $K_1$ does.
\end{proposition}
\section{Normal projective connections}
A normal projective connection is 
complicated to define precisely.
\begin{definition}
A \emph{local affine connection} on a manifold
is a choice of open set and affine connection
defined on that open set.
A \emph{covered normal projective connection} is a 
set of local affine connections whose 
open sets cover the manifold, so that 
on overlaps of the open sets, the affine connections
have the same geodesics modulo reparameterization.
A \emph{normal projective connection} is a
covered normal projective connection which is not
strictly contained in any other covered normal
projective connection.
\end{definition}
The tricky issue that makes the definition
so complicated is already seen in projective space:
projective space has no affine connection.
Each affine chart has the obvious flat affine connection. 
On the overlaps of
the affine charts, the geodesics (lines) agree.
So projective space has a normal projective connection.
\begin{definition}
The path geometry of a normal projective connection
on a manifold $S$ is the manifold $M=\Proj{}TS$,
with stalks the fibers of the obvious projection 
$\Proj{}TS \to S$, and integral curves the
curves in $\Proj{}TS$ each of which is
composed of tangent lines to a geodesic.
\end{definition}
\begin{theorem}[Hitchin \cite{Hitchin:1982}]
If the stalks [integral curves] are rational, then the path geometry
is the path geometry of a normal projective connection
on the smooth space of points [integral curves].
\end{theorem}
\begin{proof}
Fels \cite{Fels:1995} p. 239 shows that 
a 2nd order ODE system is a projective connection
(i.e. the geodesic equation of an affine
connection, in some coordinates)
just when it has the form
\begin{align*}
\frac{d^2 y}{dx^2}
&= f\left(x,y,\frac{dy}{dx}\right) \\
&= \sum_{k=0}^3 a_k(x,y) \left(\frac{dy}{dx}\right)^k.
\end{align*}
The surprising third order term arises because
$x$ is not necessarily the natural parameter
along the geodesic. Hitchin \cite{Hitchin:1982}
shows that this form
of differential equation is ensured
by rationality of the integral curves. Cartan \cite{Cartan:70}
shows that this form of differential
equation is equivalent to vanishing of $K_1$,
which turns out in local coordinates to
be a multiple of $\pd{^4 f}{\dot{y}^4}$.
\end{proof}
Theorem~\vref{thm:Main} for a scalar equation follows.
\section{Higher dimensional path geometries}
Beyond surface path geometries, the story is more complicated.
Following Mark Fels \cite{Fels:1995},
we can define a 2nd order structure on the space
of elements $M$, say $E \to M$.
We do not go through the details of the
construction of this 2nd order structure, which
is quite involved, and is explained
in detail in Fels's paper \cite{Fels:1995}.
Lets just say that the 2nd order
structure is modelled on the
tower of bundles $\PGL{n+2,\C{}} \to F \Proj{}T\Proj{n+1} \to \Proj{}T\Proj{n+1}$.
By $F \Proj{}T\Proj{n+1}$ we mean the bundle
of frames on $\Proj{}T\Proj{n+1}$, i.e.
linear isomorphism of tangent spaces 
of $\Proj{}T\Proj{n+1}$ to $\C{2n+1}$.
Tanaka \cite{Tanaka:1979} has found a Cartan connection
associated to any path geometry, using a
different normalization of torsion; this
is irrelevant to us, although it would
provide a more elegant set of
structure equations than Fels's. The interested
reader should consult \v{C}ap \cite{Cap:2005} for a clear
explanation of Tanaka's construction, 
which is superior to Fels's, but less explicit.

The structure equations of Fels's 2nd order structure are
\[
d \omega + \frac{1}{2}\left[\omega,\omega\right] = \nabla \omega,
\]
where $\omega \in \nForms{1}{E} \otimes \mathfrak{g}$
is a 1-form valued in $\mathfrak{g}=\slLie{n+2,\C{}}$, which
we can write as
\[
\omega =
\begin{pmatrix}
\omega^0_0 & \omega^0_1 & \omega^0_J \\
\omega^1_0 & \omega^1_1 & \omega^1_J \\
\omega^I_0 & \omega^I_1 & \omega^I_J \\
\end{pmatrix}
\]
where $\omega^0_0+\omega^1_1+\omega^I_I=0$,
with indices
\begin{align*}
\mu,\nu,\sigma,\tau&=0,\dots,n+1\\
i,j,k,l&=1,\dots,n+1 \\
I,J,K,L&=2,\dots,n+1 \\
\end{align*}

The terms $\nabla \omega$ have the
form
\[
\nabla \omega^{\mu}_{\nu} =
\frac{1}{2} t^{\mu}_{\nu kl} \omega^k_0 \wedge \omega^l_0
+
T^{\mu}_{\nu k l} \omega^k_1 \wedge \omega^l_0
+
\TT^{\mu K}_{\nu l} \omega^1_K \wedge \omega^l_0.
\]
These torsion terms satisfy a large collection of
identities. For example, $0=\nabla \omega^1_0=\nabla \omega^I_0$.
The torsion terms and many (perhaps all) of their identities
are worked out in Fels \cite{Fels:1995},
while Grossman \cite{Grossman:2000} assumes vanishing
of the 1-torsion of his 2nd order structure, and therefore
his identities do not cover the generality required
for our results.

The structure group of the 2nd order structure is the group
of projective transformations fixing a pointed line.
The 2nd order structure is \emph{not} necessarily a Cartan connection,
since $\nabla \omega$ is not necessarily semibasic for the
map $E \to M$, i.e. is not necessarily a multiple of the 1-forms
$\omega^1_0,\omega^I_0,\omega^I_1$ which span the
semibasic 1-forms for that map. However, the
mysterious $\TT$ terms, which precisely form the
obstruction to being a Cartan connection, vanish
in every term except $\nabla \omega^0_1, \nabla \omega^0_J$, where they
are
\begin{align*}
\TT^{0 K}_{1 L} &= \frac{1}{n-1} t^{K}_{1 L 1} \\
\TT^{0 K}_{J 1} &= \left(1+\frac{2}{n-1}\right) t^{K}_{1J1}.
\end{align*}
\begin{lemma}
Consider the 1-form $\omega$ on Fels's 2nd
order structure $E \to M$ for a path geometry
on $M$.
The obstructions to $\omega$ being a Cartan connection
vanish just precisely when the 1-torsion of the 2nd order structure vanishes,
i.e. when $\nabla \omega^1_0=\nabla \omega^I_0 = \nabla \omega^I_1=0$.
\end{lemma}
\begin{proof}
Grossman \cite{Grossman:2000} p. 435 shows that the vanishing of $t^I_{1J1}$
ensures the vanishing of $t^I_{1JK}$, and that this ensures the vanishing
of all $\TT$ terms. 
It is a long, but not difficult, calculation which requires
only differentiating the structure equations,
made much easier by using my notation.
This is precisely the condition for
vanishing of the 1-torsion of the 2nd order structure.
\end{proof}
The reader may be curious as to how to obtain
the identities. The pattern one observes in differentiating
structure equations is straightforward: we have
\[
d \omega + \omega \wedge \omega = \nabla \omega
\]
so that taking exterior derivative gives
the \emph{Bianchi identity}
\[
d \nabla \omega = \nabla \omega \wedge \omega - \omega \wedge \nabla \omega.
\]
If $\nabla \omega$ satisfies some
identity, with constant coefficients, then so
must $d \nabla \omega$, from which the above
equation gives more identities on $\nabla \omega$.
Computing out these equations gives
\begin{align*}
0 = \nabla \omega^1_0 \implies 0 &= \nabla \omega^1_I \wedge \omega^I_0 \\
0 = \nabla \omega^I_0 \implies 0 &= \nabla \omega^I_1 \wedge \omega^1_0
+ \left( \nabla \omega^I_J - \delta^I_J \nabla \omega^0_0 \right)\wedge \omega^J_0 \\
0 = \nabla \omega^0_0 + \nabla \omega^1_1 + \nabla \omega^I_I \implies 0 &= \nabla \omega^0_1 \wedge \omega^1_0
+ \nabla \omega^0_J \wedge \omega^J_0.
\end{align*}

Another observation which entails many identities
is that in the process of the method of equivalence,
covariant derivatives $\nabla$ are always semibasic.
The 1-forms $\omega^1_0,\omega^I_0,\omega^I_1$
(those under the diagonal of the matrix $\omega$)
are semibasic for $E \to M$. But since $E$ is a 2nd-order
structure, it is built on a first order structure,
say $E \to F \to M$, and the on-diagonal entries
$\omega^0_0, \omega^1_1, \omega^I_J$, together
with the below diagonal, are semibasic
for $E \to F$. Therefore
\begin{align*}
\nabla
\begin{pmatrix}
\omega^1_0 \\
\omega^I_0 \\
\omega^I_1
\end{pmatrix}
&=0
\mod \left(\omega^1_0,\omega^K_0,\omega^K_1\right)^2, \\
\nabla
\begin{pmatrix}
\omega^0_0 \\
\omega^1_1 \\
\omega^I_J
\end{pmatrix}
&=0
\mod \left(\omega^1_0,\omega^K_0,\omega^K_1,\omega^0_0,\omega^1_1,\omega^K_L\right)^2.
\end{align*}
\begin{sidewaystable}
\begin{centering}
\[
\begin{array}{rcl}
\toprule
\head{This paper} & \head{Fels} & \head{Grossman} \\
\midrule
I,J,K,L=2,\dots,n+1 & i,j,k,l=1,\dots,n & i,j,k,l=1,\dots,n \\
& i=I-1, \text{etc.} & i=I-1, \text{etc.} \\
\begin{pmatrix}
\omega^0_0 & \omega^0_1 & \omega^0_J \\
\omega^1_0 & \omega^1_1 & \omega^1_J \\
\omega^I_0 & \omega^I_1 & \omega^I_J \\
\end{pmatrix}
&
\begin{pmatrix}
\frac{1}{n+2}\left( \alpha + \Omega^i_i \right) &
\sigma &
-\beta_j \\
\omega &
\frac{1}{n+2}\left( \Omega^i_i - \left(n+1\right) \alpha \right ) &
- \kappa_j \\
\theta^i &
\pi^i &
-\Omega^i_j + \frac{\delta^i_j}{n+2}\left( \Omega^k_k + \alpha \right)
\end{pmatrix}
&
\begin{pmatrix}
\epsilon & \gamma & \pi_j \\
\lambda & \alpha + \epsilon & \zeta_j \\
\theta^i & \eta^i & \beta^i_j + \delta^i_j \epsilon
\end{pmatrix}
\\
t^I_{1J1} &
\tilde{P}^i_j &
T^i_j \\
t^I_{1JK} &
2\tilde{Q}^i_{jk} &
S^i_{jk} \\
t^0_{1K1} &
-Y_k &
\NA \\
t^0_{1KL} &
2 X_{kl} &
\NA \\
T^0_{1KL} &
W_{kl} &
\NA \\
T^0_{JKL} &
-\left(\lambda_{jk}\right)_{\theta^l} &
P_{jkl} \\
t^1_{JKL} &
-\left(\xi^2_{jk}\right)_{\theta^l} &
Z_{jkl} \\
t^{I}_{JK1} &
-\tilde{Q}^i_{jk} + T^i_{jk} - \frac{\delta^i_j}{n+2}T^l_{lk} &
\NA \\
t^I_{JKL} &
-R^i_{jkl} + \frac{\delta^i_j}{n+2} R^m_{mkl} &
\NA \\
T^I_{JKL} &
-\tilde{S}^i_{jkl} &
B^i_{jkl}  \\
\bottomrule
\end{array}
\]
\end{centering}
\caption{Dictionary of notation between this paper, Fels \cite{Fels:1995}, and Grossman \cite{Grossman:2000};
$\NA$ indicates that no notation was provided for this quantity.}
\end{sidewaystable}
\begin{proposition}
The 2nd order structure of Fels determines and is determined by a unique
path geometry. Moreover, every 2nd order structure modelled
on the tower of bundles 
$\PGL{n+2,\C{}} \to F \Proj{}T\Proj{n+1} \to \Proj{}T\Proj{n+1}$
is Fels's 2nd order structure of a path geometry if and only if
it satisfies $\nabla \omega^1_0=\nabla \omega^I_0$.
\end{proposition}
\begin{proof} We shall sketch the proof, which
depends on the details of Fels's argument in \cite{Fels:1995}.
Given the path geometry, we leave it to Fels
to construct the 2nd order structure. Given the 2nd order structure,
say $E \to M$, pick any local section $\sigma$. Following
Fels's definition of $E$ (which is complicated), one
finds that $\left(0=\sigma^* \omega^I_0=\sigma^* \omega^1_0\right)$
is the foliation by integral curves, while the foliation
by stalks is $\left(0=\sigma^* \omega^1_0 = \sigma^* \omega^I_0 \right)$.
So if there is a path geometry inducing the 2nd order structure,
then this is it. Retracing Fels's steps, we can see
that the 2nd order  structure is now completely determined,
since Fels's algorithm for constructing the 2nd order structure
depends only on having a given path geometry and
forcing the equations $0 = \nabla \omega^1_0 = \nabla \omega^I_0$,
which is enough to determine the rest of his equations
on torsion.
\end{proof}
\section{Rational integral curves}
Clearly an integral curve is rational just when it is compact with finite
fundamental group, by the classification of complex curves (see Forster \cite{Forster:1991}).
\begin{theorem}\label{thm:TF}
If the generic integral curve of a path geometry
on a manifold $M^{2n+1}$ is rational, then the path geometry is
torsion-free (i.e. its Fels/Tresse torsion vanishes).
\end{theorem}
\begin{proof}
We prove the result for $n>1$, i.e. for the Fels
torsion, since the result for $n=1$ is proven above.
Let $\iota : C \to M$ be an immersed integral curve,
$E \to M$ the bundle constructed by Fels,
with 1-form $\omega \in \nForms{1}{E} \otimes \slLie{n+2,\C{}}$.
On the pullback bundle $\iota^* E \to C$,
$\omega^I_0=0$. The structure equations simplify
greatly. Indeed $\omega$ forms a flat Cartan connection on
$C$, with $\omega \in \nForms{1}{\iota^* E} \otimes \gline$,
modelled on $\Gline/\Gptline$, the Cartan connection for a line in projective
space, which is easy to check. Taking exterior derivative
of the structure equations, we find that on $E$ the invariant $t^I_{1J1}$
satisfies
\begin{align*}
\nabla t^I_{1J1}
&=
dt^I_{1J1}
+ 2 \left( \omega^0_0 - \omega^1_1 \right) t^I_{1J1}
+ \omega^I_K t^K_{1J1} - t^I_{1K1} \omega^K_J
\\
&=
\nabla^0_K t^I_{1J1} \omega^K_0 + \nabla^1_K t^I_{1J1} \omega^K_1,
\end{align*}
for some functions $\nabla^0_K t^I_{1J1}$ and $\nabla^1_K t^I_{1J1}$.
Fixing the subalgebra $\slLie{2,\C{}} \subset \Gline$
given by the structure equations $0=\omega^0_J=\omega^1_J=\omega^I_J=\omega^0_0+\omega^1_1$,
we find
\[
dt^I_{1J1} = - 4 \omega^0_0 t^I_{1J1},
\]
so that $t^I_{1J1}$ has weight $-4$.
This ensures the vanishing of $t^I_{1J1}$ at every point of $\iota^*E.$
Since there is a raional integral curve through a generic
point of $R$, $t^I_{1J1}=0$ at all points of $E.$
Grossman \cite{Grossman:2000} p. 435 takes exterior derivatives
of the structure equations, to show
that vanishing of most of the other invariants follows,
leaving only
\begin{align*}
\nabla \omega &= d \omega + \frac{1}{2} \left[\omega,\omega\right] \\
&=
\begin{pmatrix}
0 & 0 & T^0_{JKL} \omega^K_1 \wedge \omega^L_1 \\
0 & 0 & T^1_{JKL} \omega^K_1 \wedge \omega^L_1 \\
0 & 0 & T^I_{JKL} \omega^K_1 \wedge \omega^L_1 \\
\end{pmatrix}.
\end{align*}
\end{proof}
\begin{remark}
At this stage, we may wonder if the weights
of the remaining invariants kill them as well.
Once again taking exterior derivatives,
as Grossman demonstrates,
\begin{align*}
\nabla
\begin{pmatrix}
T^0_{JKL} \\
T^1_{JKL} \\
T^I_{JKL} \\
\end{pmatrix}
=&
d
\begin{pmatrix}
T^0_{JKL} \\
T^1_{JKL} \\
T^I_{JKL} \\
\end{pmatrix}
\\
&+
\begin{pmatrix}
\left(2\omega^0_0 + \omega^1_1\right) T^0_{JKL}
- \left( T^0_{MKL} \omega^M_J + T^0_{JML} \omega^M_K + T^0_{JKM} \omega^M_L
\right)
+ \omega^0_1 T^1_{JKL}
+ \omega^0_M T^M_{JKL} \\
\left(\omega^0_0+2\omega^1_1\right) T^1_{JKL}
-\left(T^1_{MKL} \omega^M_J + T^1_{JML} \omega^M_K + T^1_{JKM} \omega^M_L \right)
\\
\left(\omega^0_0 + \omega^1_1\right)T^I_{JKL}
+ \omega^I_M T^M_{JKL}
-\left(
 T^I_{MKL} \omega^M_J + T^I_{JML} \omega^M_K + T^I_{JKM} \omega^M_L
\right)
 \\
\end{pmatrix} \\
&=
\begin{pmatrix}
0 & \nabla^0_M T^0_{JKL} & \nabla^1_M T^0_{JKL} \\
-T^0_{JKL} & \nabla^0_M T^1_{JKL} & \nabla^1_M T^1_{JKL} \\
0 & \nabla^0_M T^I_{JKL} & \nabla^1_M T^I_{JKL} \\
\end{pmatrix}
\begin{pmatrix}
\omega^1_0 \\
\omega^M_0 \\
\omega^M_1 \\
\end{pmatrix}.
\end{align*}
When I look on the same copy of $\slLie{2,\C{}}$,
structure equations turn to
\[
d
\begin{pmatrix}
T^0_{JKL} \\
T^1_{JKL} \\
T^I_{JKL} \\
\end{pmatrix}
=
\begin{pmatrix}
- \omega^0_0 T^0_{JKL}+ \omega^0_1 T^1_{JKL} \\
\omega^0_0 T^1_{JKL} - \omega^1_0 T^0_{JKL}  \\
0
\end{pmatrix}.
\]
So $T^0_{JKL}$ doesn't have a weight, since
it has a $\omega^0_1$ term, while the weights
of $T^1_{JKL},T^I_{JKL}$ are $1$ and $0$
respectively. Therefore we cannot conclude
that these invariants vanish; we soon
consider how they could come about.
\end{remark}
\subsection{Segr{\'e} geometries}
We have just seen why rationality of integral curves
forces vanishing of torsion. We need to see
why vanishing of torsion ensures rationality
of integral curves (modulo local isomorphism).
Grossman \cite{Grossman:2000}
proved that torsion-free path geometries
are locally isomorphic to path geometries
derived from 
Segr{\'e} structures,
so we need to define and examine 
Segr{\'e} structures,
to see if the derived path geometries
have rational integral curves.

In the model case, that of lines in projective
space, i.e. the system of ordinary differential equations
\[
\frac{d^2 y^I}{dx^2} = 0,
\]
the space of integral curves is
the Grassmannian of lines in projective space,
i.e. of 2-planes in $\C{n+1}$. Therefore
we should try to understand the local geometry
of the Grassmannian clearly, and look for
analogies when studying general 2nd order
systems. The Grassmannian
is $\Gtot/\Gline$ where $\Gtot=\PGL{n+1,\C{}}$,
and $\Gline$ the subgroup of $\Gtot$ fixing a projective
line.
\begin{definition}
A Cartan geometry $E \to \Dual$ modelled
on the Grassmannian of 2-planes in $\C{n+2}$
is called a \emph{Segr{\'e} geometry}.
Let $\Gpt \subset \Gtot$
be the subgroup of transformations preserving
a point on the given projective line,
and $\Gptline = \Gpt \cap \Gline$
the subgroup fixing the point and the line,
so that the model Segr{\'e} geometry is $\Gr{2}{\C{n+2}}=\Gtot/\Gline$.
The space $E/\Gptline$ is called the
\emph{space of elements} of the Segr{\'e} geometry.
The fibers of $E/\Gptline \to E/\Gline = \Dual$
are called the \emph{integral curves} of the Segr{\'e} geometry.
\end{definition}
For the Grassmannian, the space of elements is
the space of choices of a line in projective
space and a point on that line. The integral curves of the
Grassmannian are the choices of points
lying in a given line in projective space.
Keep in mind that the integral curves of a Segr{\'e} geometry
$E \to \Dual$ are \emph{not} submanifolds of the base manifold
$\Dual$, but rather the fibers of the
space of elements $M$ as a bundle $M \to \Dual$.
\begin{lemma}
The integral curves of a Segr{\'e} geometry are rational curves.
\end{lemma}
\begin{proof}
They are copies of $\Gline/\Gptline=\Proj{1}$.
\end{proof}
\begin{lemma}
The space of elements of a Segr{\'e} geometry
bears a canonical Cartan geometry modelled
on $\Gtot/\Gpt$, for which all integral curves are rational curves.
\end{lemma}
\begin{proof}
Suppose that $E \to \Dual$ is a Segr{\'e} geometry,
with Cartan connection $\omega$.
Let $M$ be the space of elements of that
Segr{\'e} geometry. One easy checks the hypotheses
of a Cartan connection to see that $\omega$ is a
Cartan connection for $E \to M$, modelled on
the space of pointed lines in projective space.
\end{proof}
\begin{definition}
Let $\Gtot' \subset \Gline$ be the subgroup
fixing a point of the Grassmannian
and fixing the tangent space to the
Grassmannian at that point.
A Segr{\'e} geometry $E \to \Dual$ is
called \emph{torsion-free} if
$\nabla \omega = 0 \pmod{\gtot'}$.
\end{definition}
Torsion-freedom of a Segr{\'e} geometry $E \to \Dual$
ensures that the equations $\omega^1_0=\omega^I_0=0$
are holonomic (i.e. satisfy the conditions of the Frobenius),
so that the manifold $E$ is foliated by the integral
manifolds of this equation. Moreover, the reader can
check that each integral manifold maps to an immersed submanifold
of $\Lambda$, called naturally a \emph{stalk} of the
Segr{\'e} geometry. Indeed the stalks foliate
the space of elements.
\begin{theorem}
A path geometry is straight just when it is
torsion-free, which occurs just when it is
locally isomorphic to the path geometry
of a unique torsion-free Segr{\'e} geometry.
\end{theorem}
\begin{proof}
We have seen that the integral curves of the
Cartan connection of the space of elements
of a Segr{\'e} geometry are
rational, hence the path geometry is straight
and therefore torsion-free. If we have a torsion-free
path geometry, then its Cartan connection
satisfies the torsion-freedom condition
required of a torsion-free Segr{\'e} geometry.
By theorem~\ref{thm:local}, the Cartan connection
is locally isomorphic to the
Cartan connection of a torsion-free Segr{\'e} geometry.
Therefore the space of elements of
the Segr{\'e} geometry is identified locally
with the path geometry.

If we have a torsion-free Segr{\'e} geometry,
then its structure equations are precisely
those of a torsion-free path geometry on the
space of pointed lines, with its integral curves as
integral curves, and stalks as stalks.
\end{proof}
Theorem~\vref{thm:Main} follows.
\begin{proposition}[Grossman \cite{Grossman:2000}]
The general torsion-free Segr{\'e} geometry
on a manifold $\Dual$ of dimension $2(n-1)$ depends
on $n(n+1)$ arbitrary functions of $n+1$ variables.
\end{proposition}
Grossman's proof unfortunately employs the Cartan--K{\"a}hler
theorem, which is not constructive. There is no
known construction producing the torsion-free
Segr{\'e} geometries, or even any large family
of examples of them. It would be very interesting
to classify the homogeneous torsion-free Segr{\'e}
geometries, and those of low cohomogeneity.
\subsection{Segr{\'e} structures}
Just as for normal projective connections,
we need to take care in defining Segr{\'e}
structures.
\begin{definition}
A \emph{local Segr{\'e} structure} on a manifold $\Dual$
of dimension $2n$ is a choice of an open set $\Omega \subset \Dual$,
two vector
bundles $U,V$ on that open set of ranks $2$ and $n$ respectively, and
an isomorphism $U \otimes V = T \Omega$.
The \emph{rank} of a tangent vector is its rank as 
a tensor in $U \otimes V$.
Two local Segr{\'e} structures are \emph{equivalent} if
they give the same ranks to all tangent vectors.
A \emph{covered Segr{\'e} structure} is a set
of local Segr{\'e} structures, equivalent
on overlaps of open sets, whose open
sets cover the manifold.
One covered Segr{\'e} structure is a \emph{refinement} of another
if it strictly contains the other.
A Segr{\'e} structure is a 
covered Segr{\'e} structure not strictly
contained in any other covered Segr{\'e}
structure.
\end{definition}
The Grassmannian of 2-planes in $\C{n+2}$
has the obvious Segr{\'e} structure, and
we can choose global vector bundles $U$ and $V$:
\[
T_P \Gr{2}{\C{n+2}} \cong P^* \otimes \left(\C{n+2}/P\right).
\]
\begin{definition}
A Segr{\'e} geometry has curvature given by
\[
\nabla \omega^{\mu}_{\nu} = K^{\mu a b}_{\nu I J}
\]
where $\mu,\nu=0,\dots,n, a,b=0,1, I,J=2,\dots,n$.
Following Machida \& Sato \cite{Machida/Sato:2000}
(who follow Tanaka \cite{Tanaka:1979}), we 
say that $\omega$ is \emph{normal} if
$K^{I01}_{0KL} + K^{I10}_{0KL} = 0$.
\end{definition}
\begin{lemma}[Machida \& Sato \cite{Machida/Sato:2000}]
A Segr{\'e} geometry determines a Segr{\'e} structure.
Conversely, a Segr{\'e} structure uniquely determines
a normal Segr{\'e} geometry, reversing the
construction. The construction of each from
the other is local and smooth.
\end{lemma}
We shall not give the proof, which is long but not conceptually
difficult, following Tanaka's interpretation of Cartan's
method of equivalence. 
The group $\SL{n+2,\C{}}$ acts in the obvious representation
on $\C{n+2}$.
The group $\Gline$ is the group of
projective transformations leaving invariant the subspace
$\C{2} \subset \C{n+2}$, and write $\C{n}$ for $\C{n+2}/\C{2}$.
Under the projection, $E \to \Lambda$, say $e \mapsto \lambda$,
$\omega$ identifies $T_m M$ with $\gtot/\gline=\C{2*} \otimes \C{n}$.
Thereby, $\omega$ determines a tensor product decomposition
on each tangent space of $\Lambda$.
One has to be a little careful, since this is not
a splitting into vector bundles defined on $\Lambda$.
Taking any local section of $E \to M$ defined
on some open subset of $M$, say $\sigma : \text{open } \subset M \to E$,
we can use this prescription to define a 
local Segr{\'e} structure on that open set.
This local prescription turns out to
determine a Segr{\'e} structure. 

The bundles $U$ and $V$ are not necessarily
globally defined, because the expression
\[
E \times_{\Gline} \C{2}
\]
doesn't make sense: $\Gline$ doesn't act on $\C{2}$,
being only a subgroup of $\Gtot=\PSL{n+2,\C{}}$.
We can define the
fiber bundles $E \times_{\Gline} \Proj{1}$
and $E \times_{\Gline} \Proj{n-1},$
which we think of intuitively as
$\Proj{}U$ and $\Proj{}V$.

The distinction between local and global Segr{\'e}
geometries is not always clearly made, nor is
the distinction between Segr{\'e} structures
and Segr{\'e} geometries. The space of elements of a Segr{\'e}
structure is the total space of the bundle
$E \times_{\Gline} \Proj{1} \to \Dual$, and the fibers are the integral curves.
\begin{proposition}[Grossman \cite{Grossman:2000} p. 415]
Given a Segr{\'e} geometry $E \to \Dual$,
the 1-forms $\omega^I_0,\omega^I_1$ are
semibasic. The symmetric 2-tensors
\[
\Delta^{IJ} = \omega^I_0 \omega^J_1 + \omega^J_1 \omega^I_0
-  \omega^J_0 \omega^I_1 - \omega^I_1 \omega^J_0
\]
descend from each point of $E$ to determine
symmetric 2-tensors at the corresponding point
of $\Dual$. Their span is independent of the
choice of point in $E$, depending only on the
corresponding point of $\Dual$, defining a smooth
vector subbundle of $\Sym{2}{T \Dual}$.
\end{proposition}
\begin{proof}
It is an easy calculation that the $\Delta^{IJ}$
transform under the action of $\Gline$ as combinations
of one another just when the torsion vanishes, since
we know how $\omega$ transforms by definition
of a Cartan connection; for details see Grossman
\cite{Grossman:2000} p. 416.
\end{proof}
\begin{definition}
Let $U^2$ and $V^n$ be vector spaces of dimensions $2$ and $n$ respectively.
The \emph{Segr{\'e} variety} $\Sigma \subset \Proj{}\left(U \otimes V\right)$
is the set of elements of rank 1, i.e. pure tensors $u \otimes v$. In coordinates
$u_0,u_1$ on $U$, and $v^I$ on $V$, we have coordinates $w_0^I,w_1^I$ on
$U \otimes V$ and the Segr{\'e} variety is cut out by the equations
$w^I_0 w^J_1 = w^J_0 w^I_1$.
\end{definition}
\begin{definition}
The group $G(U,V) = \left(\GL{U} \times \GL{V}\right)/\Delta$
(where $\Delta$ is the group of
pairs of scalar multiples of the identity
of the form $\left(\lambda,\lambda^{-1}\right)$)
clearly acts as linear transformations
on $U \otimes V$ leaving the
Segr{\'e} variety invariant. If $\dim U = \dim V$,
we can also take
any linear isomorphisms $\phi,\psi : U \to V$,
and map $g(u \otimes v)=\psi^{-1}(v) \otimes \phi(u)$,
giving an action of 
\[
G'(U,V) = G(U,V) \sqcup \left(\Iso{U}{V} \times \Iso{U}{V}\right)/ \Delta.
\]
\end{definition}
\begin{lemma}
The group of linear transformations of 
$U \otimes V$ preserving the Segr{\'e} variety
is $G(U,V)$, unless $\dim U = \dim V$, in which case it is
$G'(U,V)$.
\end{lemma}
The proof is just some linear algebra.
\begin{corollary}
A Segr{\'e} structure on an even dimensional
manifold (not of dimension 4) is equivalent to a choice of
a smoothly varying family of subvarieties
in the projectivized tangent spaces of the manifold,
with each variety linearly isomorphic to the Segr{\'e} variety.
Equivalently, a Segr{\'e} structure is equivalent to a 
choice of a smoothly varying linear subspace
of the symmetric 2-tensors isomorphic at each point to the
subspace spanned by the equations cutting
out the Segr{\'e} variety.
\end{corollary}
For 4-dimensional manifolds, analoguous to
$\Gr{2}{\C{4}}$, we can consider a Segr{\'e} structure
to be a choice of a family of Segr{\'e} varieties
in the projectived tangent spaces, together with
an analogue of an orientation, picking out which
of the two tensor product factors is which.

Grossman took this view of Segr{\'e} geometries,
as families of Segr{\'e} varieties, which seems
quite natural. Nonetheless, it is not clear which
point of view makes easier the process of
geometrically constructing all of the torsion-free
Segr{\'e} geometries, a task which has yet to be done.
\begin{proposition}[Machida \& Sato \cite{Machida/Sato:2000}]
Every Segr{\'e} structure determines and is determined
by a unique normal Segr{\'e} geometry, through a
local construction. In particular, the concept of
torsion-freedom is defined for Segr{\'e} structures.
\end{proposition}
\section{Integrability}
Generic straight 2nd order ODE systems are integrable
by geometric construction, as we shall see (also see Grossmann \cite{Grossman:2000}).
Torsion yields an explicit test for integrability;
for example, every straight 2nd order ODE (i.e. $n=1$)
\[
\frac{d^2 y}{dx^2} = f\left(x,y,\frac{dy}{dx}\right)
\]
for which $f$ is a sum of linear and quadratic terms in $x,y,\dot{y}$
is integrable by use of hypergeometric functions and quadratures.
Indeed, Cartan \cite{Cartan:1938}
can apparently integrate any straight 2nd order ODE
by differentiation and at most two quadratures. (But see 
subsection~\vref{subsec:Sour}
for some concerns about Cartan's claim.)
Straight systems remain straight under symmetry reduction
so they form a fascinating class of ODE systems.

The ability to integrate straight equations is particularly
exciting when we realize that all of the 2nd order ODEs
of classical mathematical physics are straight (see
table~\vref{table:ClassicalEquations} and thus
apparently solved by Cartan's method. Special function
theory is just one special case of the theory
of straight equations. There are
some well known modern equations of mathematical
physics which are not integrable without adjoining new
transcendental functions, even if we allow inversions of
integrals of elementary functions, and thus are not obtained
by algebra and quadratures (see table~\vref{table:notStraight}); therefore they are not
integrable by Cartan's method, or by symmetry reduction.
Consonant with Cartan's claim, every such equation so far
tested has torsion. Indeed from the table, we see that torsion is
found except in certain special cases. These
cases turn out to all be well known degeneracies in which
the general solution can be expressed in hypergeometric or elliptic functions.

If Cartan's claim 
is correct, then his methods integrate in quadratures every 
integrable 2nd order ODE known to me.
Cartan's method appears to be the solution to the problem
of integration in quadratures for a single 2nd order ODE.
Geometry (rationality of integral curves) yields integrability.


\newcommand{\pp}{\frac{d^2 y}{dx^2}}
\newcommand{\p}{\frac{dy}{dx}}
\begin{table}
\begin{longtable}{p{2in}p{3in}} 
\toprule
\head{Common name} & \head{Equation} \\
\midrule
  \endfirsthead
\toprule
\head{Common name} & \head{Equation} \\
\midrule
  \endhead
  Airy &
  $\pp = xy$ \\
  Anger &
  $\pp + \frac{\p}{x} + \left(1-\frac{a^2}{x^2}\right)y = \frac{x - a}{\pi x^2} \sin \pi a$ \\
  Bessel &
  $x^2 \pp + x \p + (x^2 + a)y = 0 $ \\
  Bessel (modified)&
  $x^2 \pp + x \p - (x^2 + a)y = 0 $ \\
  Bessel (spherical) &
  $x^2 \pp + 2 x \p +(x^2 + a)y = 0 $ \\
  Bessel (modified spherical) &
  $x^2 \pp + 2 x \p - (x^2 + a)y = 0 $ \\
  confluent hypergeometric &
  $x \pp + \left(c - x\right) \p - ay = 0$ \\
  Coulomb wave &
  $\pp + \left(1-\frac{a}{x} - \frac{b}{x^2}\right)y=0 $ \\
  Eckart &
  $\pp + \left(\frac{a \, e^{dx}}{1+e^{dx}}  + \frac{b \, e^{dx}}{\left(1+e^{dx}\right)^2} + c\right)y = 0$ \\
  ellipsoidal &
  $\pp = \left(a+b\sin(x)^2+c\sin(x)^4\right)y$
  \\
  elliptic &
  $x(1-x^2) \pp + \left(1-x^2\right) \p - 2 x^2 \p - xy$
  \\
  error function &
  $\pp + 2x \p =  2ay$ \\
  Euler &
  $x^2 \pp + a x \p + b y = 0$ \\
  Gau{\ss} hypergeometric &
  $x(x-1)\pp + \left( \left(a+b+1\right)x - c \right) \p + aby = 0$ \\
  Halm &
  $\left(1+x^2\right)^2+ \pp + a \p = 0$ \\
  Hermite &
  $\pp + 2x \p = 2ay$ \\
  Lienard &
  $\pp + (ay+b) \p + \left(\frac{1}{9}a^2 y^3 + \frac{1}{3} aby^2 + cy + d\right) = 0$ \\
  Liouville &
  $\pp + g(y) \left(\p\right)^2 + f(x) \p = 0$ \\
  Mathieu &
  $\pp + \left(a-2b \cos 2x \right) y = 0$ \\
  Titchmarsh &
  $\pp + \left(b - x^a \right) y = 0$ \\
\bottomrule
\end{longtable}
\caption{Some straight equations from mathematical physics; $a,b,c,d$ any constants, $f,g$ any functions. See
Polyanin \& Zaitsev \cite{Polyanin/Zaitsev:2003}, Zwillinger \cite{Zwillinger:1992}}\label{table:ClassicalEquations}
\end{table}
\begin{table}
\begin{alignat*}{3}
\toprule
\text{\head{Common name}} & & & \text{\head{Equation}} & & \text{\raggedright{\head{When torsion is found}}} \\
\midrule
  \text{Emden--Fowler} & &
  \quad  \pp & + f(x) \p + y^a
  & & \quad a  \ne 0,1 \\
  \text{Lagerstrom} & &
  \quad x \pp & + a \p + b x y \p = 0
  & & \quad b \ne 0 \\
  \text{Painlev{\'e} I} & &
  \quad \pp &= 6y^2+x
  & & \\
  \text{Painlev{\'e} II} & &
  \quad \pp &= 2y^3+xy+a
  & & \\
  \text{Painlev{\'e} III} & &
  \quad \pp &= \frac{\p ^2}{y}-\frac{\p}{x} + \frac{ay^2+b}{x}+cy^3+\frac{d}{y}
  & \quad & (a,b,c,d) \ne (0,0,0,0) \\
  \text{Painlev{\'e} IV} & &
  \quad \pp &= \frac{\p^2}{2y} + \frac{3y^2}{2} + \\
   & & & 4y^3 x + 2\left(x^2-a\right)y + \frac{b}{y}
  & & \\
  \text{Painlev{\'e} V} & &
  \quad \pp &= \left(\frac{1}{2y} + \frac{1}{y-1} \right) \p^2
  \\
  & & &
  - \frac{\p}{x} + \frac{\left(y-1\right)^2 \left(ay+\frac{b}{y}\right)}{x^2} \\
  & & &
  + \frac{cy}{x} + \frac{dy(y+1)}{y-1}
  & \quad &(a,b,c,d)  \ne (0,0,0,0)\\
  \text{Painlev{\'e} VI} & &
  \quad \pp &= \frac{1}{2}\left( \frac{1}{y} + \frac{1}{y-1} + \frac{1}{y-x} \right) \p^2
 & &(a,b,c,d)\ne\left(0,0,0,\frac{1}{2}\right)
  \\
  & & &
  - \left( \frac{1}{x} + \frac{1}{x-1} + \frac{1}{y-x} \right) \p
  \\
  & & &
  + \frac{y(y-1)(y-x)\Gamma}{x^2(x-1)^2} \\
  & &
  \Gamma &=  a+\frac{bx}{y^2} + \frac{c(x-1)}{\left(y-1\right)^2} + \frac{dx(x-1)}{(y-x)^2}
 \\
 \text{van der Pol} & &
 \quad \pp &= a \left(1 - y^2\right) \p - y
 & & \quad a  \ne 0 \\
\bottomrule
\end{alignat*}
\caption{Some equations of mathematical physics which are not straight; $a,b,c,d$ any constants}\label{table:notStraight}
\end{table}
\subsection{A note of sour skepticism}
\label{subsec:Sour}
If a 2nd order ODE
has a Lie group of symmetries of positive dimension, it
would appear to invalidate Cartan's approach (which we shall
see in section~\vref{subsec:Integrability}) as Cartan
describes it, since the
invariants do not provide enough conservation laws. Cartan
does not point out this case, but integrability still follows
as long as the Lie group
has dimension 2 or greater (see Lie \cite{Lie:1883}).
Even if the symmetry group is not solvable, the equation is
integrable. For example, consider
\[
\frac{d^2 y}{dx^2} = 0,
\]
whose symmetry Lie algebra turns out to be $\slLie{3,\C{}}$,
a simple Lie group. However, one needs to know the symmetry
group action explicitly to carry out this integration.
The question of the integrability of a
straight 2nd order ODE in the presence of a one dimensional
Lie group of symmetries is apparently not settled, in
contrast to Cartan's claim.

The equations of mathematical physics described
above, as a consequence of the theorems we have proven, 
are all locally equivalent under point transformations
to the standard equation $d^2 y/dx^2 = 0$, and therefore
have simple Lie pseudogroup of point symmetries, so that
Lie's method of reduction does not apply. The symmetry
groups are not explicit, and finding them explicitly
appears to be as difficult as solving the equations
directly. Moreover,
Cartan's approach as he outlines it also does not apply,
since it depends on local invariants under point transformations.
It may be that Cartan has up his sleeve some deeper methods
that apply in these ``degenerate'' cases, but he gives
no indication. Nonetheless, it is amazing that
the basic ODEs of mathematical 
physics (before Painlev\'e) are straight,
that straightness is a rare property,
and that no one has previously noticed this.
\subsection{Optimism returns, with topology in tow}
The theory of 2nd order ODEs of mathematical
physics appears from this point of view to be
nearly topological, in the sense that all of the
straight equations of mathematical phyiscs
are locally point equivalent to
$d^2 y/dx^2 =0$, i.e. to the contact 3-manifold
$x,y,p$ with contact planes $dy = p \, dx$
and two Legendre foliations (a) $dy = p \, dx, dp = 0$
and (b) $dx=dy=0$. The global study of such ``flat''
double Legendre folations is thus at the core of
physics, while being locally completely elementary.
Contact topology with flat
double Legendre folation is entirely mysterious.
\subsection{Why we study 2nd order systems, not first order ones}
All first order systems of ordinary differential equations
\[
\frac{dy^I}{dx} = f^I(x,y)
\]
are straight, since the Frobenius theorem tells
us that we can change coordinates to arrange
that $f^I=0$. Moreover, higher order systems can
be rewritten as first order systems, so
it might appear that they are always straight.
But this is not the case, since point
transformations of 2nd order
equations are not quite so powerful.

\subsection{Grossman's results on integrability}
\label{subsec:Integrability}
Grossman \cite{Grossman:2000} considered in some detail
the question of integrability for torsion-free path geometries.
We shall summarize his results, which
generalize Cartan's \cite{Cartan:1938}. Each torsion-free path geometry
comes from a torsion-free Segr{\'e} structure.
This Segr{\'e} structure has an associated normal Segr{\'e} geometry. 
Fels \cite{Fels:1995} shows us how to compute the structure equations
of the 2nd order structure $E \to M$ , which are precisely the
structure equations of the Segr{\'e} geometry $E \to \Lambda$. Therefore,
even though we don't see how to construct explicitly
the base manifold $\Lambda$ of the Segr{\'e} structure,
i.e. the space of solutions, we can compute its curvature, which lives on $E$.
Just by differentiating, we can compute the covariant
derivatives of all orders of the curvature. Each of
these invariants transforms under the structure group
$\Gline$ of $E \to \Lambda$  in some representation.
If we can cobble together a rational invariant
(out of these covariant derivatives) which lives in the
trivial representation of $\Gline$, then the invariant
descends to a function on the unknown manifold
$\Lambda$, i.e. on the space of integral curves, and therefore
it must be a constant on each integral curve. As Cartan
and Grossman prove, this process generically
succeeds, because there are rational invariants arising in this
manner which, for generic torsion-free Segr{\'e} structures,
have differential nonzero at a generic point. Indeed,
in this manner one can find enough conservation laws to reduce
the determination of integral curves to quadrature, integrating
the original system of ordinary differential equations. Thus
we have ``integrated by differentiating.'' This process can fail,
but only when too many invariants of the curvature and its
covariant derivatives are constant on $E$. For scalar equations
(i.e. surface path geometries), Cartan's methods 
\cite{Bryant/Griffiths/Hsu:1995,Cartan:161}
show that every torsion-free path geometry 
on a 3-manifold either has a conserved quantity, or 
the differential equation has a positive dimensional Lie group of symmetries, 
so we can reduce the equation using
Lie's method. More complicated phenomena are observed
in Grossman's thesis \cite{Grossman:2000}, where the constancy
of one particular invariant, at least in low dimension, allows
one to calculate further higher order invariants which
generically still ensure integrability. However, in general it is unknown
whether every torsion-free system of equations must
either be integrable with differential invariants as
conservation laws or have a positive dimensional
Lie group of symmetries.

\section{Rational stalks}
Given a path geometry on a complex
manifold $M^{2n+1}$, let $\Sigma \subset M$ be a stalk.
Take the Cartan geometry $E|_{\Sigma} \to \Sigma$,
which is modelled on $\Gpt/\Gptline$.
The $\omega^I_1$ are semibasic for this bundle, while
$\omega^1_0=\omega^I_0=0$. But at least one $\omega^1_0$ or $\omega^I_0$
term appears in all of the curvature of $E \to M$. Therefore
$E|_{\Sigma}$ is flat.
\begin{theorem}
A stalk of a path geometry on a complex manifold $M^{2n+1}$ is rational
just when it is compact with fundamental group defying $\Gpt$, and this occurs just
when its Cartan geometry is isomorphic to the Cartan geometry of the stalks of the model.
\end{theorem}
\begin{proof}
The Cartan connection is flat, so 
by theorem~\vref{thm:Defiance}
our compact stalk must be a locally Klein geometry
$\Gamma \backslash \Gpt / \Gptline$. But
$\Gpt / \Gptline = \Proj{n-1}$, so $\Gamma$ must be
a discrete group of projective linear transformations
acting as deck transformations on projective space.
However, every linear transformation has an eigenspace,
so every projective linear transformation has a fixed
point. Therefore $\Gamma=\left\{1\right\}$.
\end{proof}
Theorem~\vref{theorem:RatStalk} follows.
\begin{lemma}
The normal bundle of a rational stalk (as a submanifold of $M$) is trivial
$\OOp{0}{n}$.
\end{lemma}
\begin{proof}
First consider the case of the model. Each $\Proj{n}$ fiber of $\Proj{}T\Proj{n+1}$ lives
inside the open set $\Proj{}T \mathbb{A}^{n+1} = \mathbb{A}^{n+1}
\times \Proj{n}$, so clearly has trivial normal bundle
$\nu \Proj{n} = \OOp{0}{n+1}$. Next, in the general
case, construct the normal bundle as
\[
\Norm{\Sigma}{M} = \left(E|_{\Sigma} \times \left( \gpt/\gptline \right)\right)/\Gptline.
\]
Therefore $\Norm{\Sigma}{M} = \OOp{0}{n+1}$. But $E|_{\Sigma} \to \Sigma$ is
isomorphic to the model $\Gpt \to \Proj{n}$.
\end{proof}
\begin{theorem}
If the space of elements of a path geometry is connected,
and all stalks are compact, and one stalk has fundamental group
defying $\Gpt$, then all stalks are rational,
and the space of points is a smooth complex manifold, and
the map taking an element to its point is smooth.
\end{theorem}
\begin{proof}
Follows immediately from Kodaira theory.
\end{proof}
\begin{lemma}
If the stalks of a path geometry are rational, then the invariant $T^I_{JKL}$
vanishes.
\end{lemma}
\begin{proof}
Following Fels \cite{Fels:1995} p. 235, we compute
that
\[
\nabla T^I_{JKL}
=
dT^I_{JKL} +
\left( \omega^0_0 + \omega^1_1 \right) T^I_{JKL}
+
\omega^I_M T^M_{JKL}
- T^I_{MKL} \omega^M_J
- T^I_{JML} \omega^M_K
- T^I_{JKM} \omega^M_L
\]
is semibasic for the map $E \to M$. Pick a
number $N$ from $2,\dots,n$.  Consider the copy of $\slLie{2,\C{}} \subset \Gpt$
given by the equations $\omega^1_1 + \omega^N_N=0$
together with setting every $\omega^{\bullet}_{\bullet}$
to $0$ except for
$\omega^{1}_{1},\omega^{N}_{1},\omega^{1}_{N},\omega^{N}_{N}$.
Calculate that
\[
dT^I_{JKL} = T^I_{JKL} \omega^1_1 \left(\delta^I_N - 1 - \delta^N_J - \delta^N_K - \delta^N_L \right).
\]
If $I \ne N$, then clearly this is a negative line bundle.
Therefore $T^I_{JKL}=0$ as long as $I \ne N$. But if
$I = N$, then switch to a different choice of index $N$.
\end{proof}
\begin{theorem}
All of the stalks of a path geometry
are rational just when then the space of
points is a smooth complex manifold, and the
integral curves of the path geometry project to the
geodesics of a unique normal projective connection.
In particular, near any point of the point space,
the projected integral curves are precisely the
geodesics of some affine connection.
\end{theorem}
\begin{proof}
Kodaira's theorem ensures that the space of
stalks is a complex manifold.
Following Fels \cite{Fels:1995} p. 238,
the vanishing of $T^I_{JKL}$ is precisely
the condition under which the path geometry
is locally that of a projective connection
on some complex manifold, which is locally
identified with Kodaira's moduli space.
\end{proof}
We now prove theorem~\vref{theorem:Local}.
\begin{proof}
This is a long calculation: once the invariants
$T^I_{JKL}$ and $\TT$ are forced to vanish, then
all of the remaining invariants vanish, and then
the Cartan geometry on $E$ is flat.
\end{proof}
We now prove theorem~\vref{thm:Model}.
\begin{proof}
The stalks are rational, so the path geometry
is a normal projective connection on the space
of points $S$, which is a smooth manifold.
The normal projective connection is flat,
so a covering space $\tilde{S}$ of $S$ 
is mapped to projective space, and
the normal projective connection pulled back. 
The integral curves of the path geometry project
to the geodesics, so the geodesics are rational
curves.
Because each geodesic is simply connected,
and admits no smooth quotient curve,
each geodesic in $\tilde{S}$ maps bijectively to a
projective line. Since any two points in projective space 
lie on a projective line, the map
$\tilde{S} \to \Proj{n+1}$ is a surjective
local diffeomorphism.

Take a point $s \in \tilde{S}$, and suppose
it is mapped to a point $p \in \Proj{n+1}$.
Let $B$ be the blowup of $\Proj{n+1}$ at $p$.
So points of $B$ are 
pairs $\left(\ell,q\right)$ with $\ell$
a line through $p$ and $q$ a point of that
line. Given $(\ell,q)$, let $\tilde{\ell}$ be 
the geodesic through $s$ which is mapped
to $\ell$. Since the map $\tilde{S} \to \Proj{n+1}$
is bijective
on each geodesic, there is a unique
point $\tilde{q}$ on $\tilde{\ell}$
mapping to $q$. Map $(\ell,q) \in B \to \tilde{q} \in \tilde{S}$.
Clearly the map has image consisting precisely
in the points which lie on a geodesic through $s$.
Moreover, $B$ is compact, so the image
of this map must be as well. Therefore the
geodesics through any point of $\tilde{S}$
cover a compact subset of $\tilde{S}.$
The composition $B \to \tilde{S} \to \Proj{n+1}$ is 
the blowup map, so a local biholomorphism
on a dense open set. Therefore $B \to \tilde{S}$
is holomorphic, and on some open set a local biholomorphism.
By Sard's lemma, the map $B \to \tilde{S}$
is onto. Therefore $\tilde{S}$ is compact,
and the map $\tilde{S} \to \Proj{n+1}$
is a biholomorphism.
\end{proof}

\section{Literature}
The approach we take here is very similar
to Hitchin \cite{Hitchin:1982} and Dunajski \& Todd
\cite{DunajskiTodd:2005}. Merker \cite{Merker:2004} has a different approach,
which characterizes very explicitly the
systems of ordinary differential equations
equivalent to $d^2y/dx^2=0$. 
The papers of Fritelli, Simonetta, Kozameh and Newman 
\cite{Frittelli/Simonetta/Kozameh/Newman:2001}
and of Newman and Nurowski \cite{Newman/Nurowski:2003}
concern questions closely related to this paper.
Bordag \& Dryuma \cite{BordagDryuma:1997} make use of Cartan's invariants
of projective connections to analyse 2nd order ordinary
differential equations.

\section{Open problems}
Some open problems:
\begin{enumerate}
\item
Write software to symbolically integrate ``generic'' torsion-free
systems of ordinary differential equations.
\item
For torsion-free systems for which there are not enough
invariants to integrate (using curvature and its covariant
derivatives to generate integrals of motion), is there always
some other process to integrate the equations, by combination
of those integrals of motion together with symmetry reductions?
\item
Find a straightness criterion for higher order systems of equations,
for example third order scalar equations \cite{Cartan:174,Chern:1940,Sato/Yoshikawa:1998},
fourth order scalar equations \cite{Bryant:1991,Fels:1996}, and third order systems
\cite{Fels:1993}. Dunajski and Todd \cite{DunajskiTodd:2005} have recent solved
this problem for any $n$-th order ODE.
\item
One can adapt the methods of this article to a host of
Cartan connections and $G$-structures. For example to 2-plane
fields on a 5-manifold, satisfying a natural nondegeneracy
condition (see Cartan \cite{Cartan:30}): one can ask when their
bicharacteristic curves are rational. The crucial idea
is to look at a copy of $\slLie{2,\C{}}$ appearing in the
structure equations, and see how the torsion (or curvature)
varies under it, which we can read off directly from
structure equations.
\item
Can something be said about ordinary differential equations
whose integral curves are elliptic? The methods employed
here seem powerless, since line bundles on elliptic
curves have moduli, so we couldn't expect to read
them off from the structure equations.
\item
The requirement that a Segr{\'e} structure be torsion-free
is a collection of first order partial differential equations,
which has a lot of local solutions (the Cartan--K{\"a}hler
theorem tells us so). But there is no technique for
constructing solutions. We are not interested in flat
solutions (i.e. locally isomorphic to the Grassmannian
of 2-planes in a vector space), but quite interested
to find the nonflat examples with largest possible symmetry
groups, which correspond to very special systems
of ordinary differential equations.
\item
Cartan's concept of ``integrating by differentiating''
applies to certain families of ordinary differential
equations, which he refered to as \emph{classe C}
\cite{Cartan:1938}. Is there actually a relation
between straightness and class C? Presumably
straightness implies class C, but comments
in Bryant \cite{Bryant:1991} p. 35 suggest that there
class C might not imply straightness. 
\item
In a subsequent paper, I will
demonstrate constraints on the 
characteristic classes of
closed K\"ahler manifolds
admitting path geometries.
In another paper, I will classify the projective 3-folds which
admit path geometries, just as
Jahnke and Radloff \cite{Jahnke/Radloff:2002,Jahnke/Radloff:2004} did
for normal projective connections
and conformal structures.
\item
Perhaps if the path geometry is singular,
but all of the integral curves are still rational,
there is still some local information.
\end{enumerate}
\bibliographystyle{amsplain}
\bibliography{Grossman}
\end{document}